%% file: RareOpt.tex
\documentclass[final,onefignum,onetabnum]{siamonline190516}

\input{Preamble.tex}

\begin{document}
\maketitle

\begin{abstract}
Chance constraints provide a principled framework to mitigate the risk of high-impact extreme events by modifying the controllable properties of a system. The low probability and rare occurrence of such events, however, impose severe sampling and computational requirements on classical solution methods that render them impractical. This work proposes a novel sampling-free method for solving rare chance constrained optimization problems affected by uncertainties that follow general Gaussian mixture distributions. By integrating modern developments in large deviation theory with tools from convex analysis and bilevel optimization, we propose tractable formulations that can be solved by off-the-shelf solvers. Our formulations enjoy several advantages compared to classical methods: their size and complexity is independent of event rarity, they do not require linearity or convexity assumptions on system constraints, and under easily verifiable conditions, serve as safe conservative approximations or asymptotically exact reformulations of the true problem.  Computational experiments on linear, nonlinear and PDE-constrained problems from applications in portfolio management, structural engineering and fluid dynamics illustrate the broad applicability of our method and its advantages over classical sampling-based approaches in terms of both accuracy and efficiency.
\end{abstract}



\input{introduction.tex}

\input{literature.tex}

\input{large_deviation.tex}

\input{gaussian.tex}

\input{gaussian_mixture.tex}

\input{applications.tex}

\input{conclusions.tex}

\section*{Acknowledgments}
We thank Eric Vanden-Eijnden and Georg Stadler from Courant Institute at New York University for helpful suggestions at various stages of this work.

\bibliographystyle{siamplain}
\bibliography{RareOpt,References}

\end{document}

%% file: Preamble.tex
\usepackage{braket,amsfonts}

\usepackage{array}

\usepackage[caption=false]{subfig}
\usepackage{pgfplots}

\usepackage{algorithmic}

\usepackage{graphicx,epstopdf}

\usepackage{tikz,pgfplots}
\usetikzlibrary{spy}
\usepackage{makecell}
\usepgfplotslibrary{fillbetween}
\usetikzlibrary{pgfplots.groupplots}
\usepackage{multirow}

\usepackage{mathtools}

\usepackage{cleveref}
\usepackage{dsfont}
\usepackage{esint}
\usepackage{color,cite}
\usepackage{amssymb}
\usepackage{enumitem}

\usepackage{booktabs,tabularx}
\allowdisplaybreaks

\Crefname{ALC@unique}{Line}{Lines}

\usepackage{amsopn}

\usepackage{xspace}
\usepackage{bold-extra}
\usepackage[most]{tcolorbox}

\colorlet{texcscolor}{blue!50!black}
\colorlet{texemcolor}{red!70!black}
\colorlet{texpreamble}{red!70!black}
\colorlet{codebackground}{black!25!white!25}

\lstdefinestyle{siamlatex}{%
  style=tcblatex,
  texcsstyle=*\color{texcscolor},
  texcsstyle=[2]\color{texemcolor},
  keywordstyle=[2]\color{texemcolor},
  moretexcs={cref,Cref,maketitle,mathcal,text,headers,email,url},
}

\tcbset{%
  colframe=black!75!white!75,
  coltitle=white,
  colback=codebackground, %
  colbacklower=white, %
  fonttitle=\bfseries,
  arc=0pt,outer arc=0pt,
  top=1pt,bottom=1pt,left=1mm,right=1mm,middle=1mm,boxsep=1mm,
  leftrule=0.3mm,rightrule=0.3mm,toprule=0.3mm,bottomrule=0.3mm,
  listing options={style=siamlatex}
}

\newtcblisting[use counter=example]{example}[2][]{%
  title={Example~\thetcbcounter: #2},#1}

\newtcbinputlisting[use counter=example]{\examplefile}[3][]{%
  title={Example~\thetcbcounter: #2},listing file={#3},#1}

\DeclareTotalTCBox{\code}{ v O{} }
{ %
  fontupper=\ttfamily\color{black},
  nobeforeafter,
  tcbox raise base,
  colback=codebackground,colframe=white,
  top=0pt,bottom=0pt,left=0mm,right=0mm,
  leftrule=0pt,rightrule=0pt,toprule=0mm,bottomrule=0mm,
  boxsep=0.5mm,
  #2}{#1}

\patchcmd\newpage{\vfil}{}{}{}
\def\EE{\mathbb{E}}
\def\PP{\mathbb{P}}
\def\RR{\mathbb{R}}

\DeclareRobustCommand{\argmin}{\operatorname*{argmin}}

\DeclareRobustCommand{\argmax}{\operatorname*{argmax}}

\DeclareRobustCommand{\CVaR}{\operatorname*{CVaR}}
\DeclareRobustCommand{\VaR}{\operatorname*{VaR}}

\def\EE{\mathbb{E}}\def\PP{\mathbb{P}}
\def\RR{\mathbb{R}}

\def\<{\left\langle} \def\>{\right\rangle}

\newcommand{\Gauss}[2]{\mathcal{N}(#1,#2)}
\newcommand{\xistar}{\xi^{\star}}

\newcommand{\etastar}{\eta^{\star}}

\newcommand{\SigInv}[1]{\Sigma_{#1}^{-1}}
\newcommand{\SigInvHalf}[1]{\Sigma_{#1}^{-\frac{1}{2}}}
\newcommand{\SigHalf}[1]{\Sigma_{#1}^{\frac{1}{2}}}
\newcommand{\nhat}{\hat{n}}

\newcommand{\indicator}[1]{\mathds{1}_{#1}}
\newcommand{\nstar}{\hat{n}^\star}
\newcommand{\Pm}{\mathcal{P}}
\newcommand{\txi}{\tilde{\xi}}
\newcommand{\tn}{\tilde{n}}

\newcommand{\lwsc}{1pt} %
\newcommand{\mssc}{2pt} %

\newcommand{\mLDTF}{square*} %
\newcommand{\mLDTS}{*} %
\newcommand{\mSAT}{triangle*} %
\newcommand{\mCVaR}{diamond*} %

\newcommand{\cLDTF}{orange!80!white} %
\newcommand{\cLDTS}{red!80!black} %
\newcommand{\cSATtwo}{blue} %
\newcommand{\cSATthree}{violet} %
\newcommand{\cSATfour}{gray} %
\newcommand{\cSATfive}{cyan} %
\newcommand{\cCVaRfive}{green!70!black} %

\pgfplotscreateplotcyclelist{my color}{%
line width=\lwsc,mark size=\mssc,	solid, \cLDTF!50!black,every mark/.append style={solid, fill=\cLDTF}, mark=\mLDTF \\ %
line width=\lwsc,mark size=\mssc,	densely dashed,\cLDTS!50!black,every mark/.append style={solid, fill=\cLDTS}, mark=\mLDTS\\ %
line width=\lwsc,mark size=\mssc,	solid, \cCVaRfive!50!black,every mark/.append style={solid, fill=\cCVaRfive}, mark=\mCVaR\\ %
line width=\lwsc,mark size=\mssc,	dotted, \cSATtwo!50!black, every mark/.append style={solid, fill=\cSATtwo}, mark=\mSAT\\
line width=\lwsc,mark size=\mssc,	dashed, \cSATthree!50!black,every mark/.append style={solid, fill=\cSATthree},mark=\mSAT\\
line width=\lwsc,mark size=\mssc,	dashdotted, \cSATfour!50!black, every mark/.append style={solid, fill=\cSATfour},mark=\mSAT\\
line width=\lwsc,mark size=\mssc,	solid, \cSATfive!50!black,every mark/.append style={solid, fill=\cSATfive},mark=\mSAT\\
	dashdotted, every mark/.append style={solid, fill=gray},mark=otimes*\\%
	dashdotdotted, every mark/.append style={solid},mark=star\\%
	densely dashdotted,every mark/.append style={solid, fill=gray},mark=diamond*\\%
}

\pgfplotscreateplotcyclelist{my mark}{
	green!30!black,every mark/.append style={fill=green!50!black},mark=*\\
	teal!60!black,every mark/.append style={fill=teal!80!black},mark=square*\\ 
	cyan!30!black,every mark/.append style={fill=cyan!50!black},mark=diamond*\\ 
	cyan!60!black,every mark/.append style={fill=cyan!80!blue},mark=pentagon*\\
	blue!30!black,every mark/.append style={fill=blue!50!white},mark=otimes*\\
 violet!60!black,	every mark/.append style={solid,fill=violet!80!black},mark=triangle*\\ yellow!60!black,densely dashed,
	every mark/.append style={solid,fill=yellow!80!black},mark=halfsquare*\\
	black,every mark/.append style={solid,fill=gray},mark=halfcircle*\\
	blue,densely dashed,mark=star,every mark/.append style=solid\\
	red,densely dashed,every mark/.append style={solid,fill=red!80!black},mark=halfdiamond*\\
}

\newcommand{\edits}[1]{#1}
\newcommand{\modifications}[1]{#1}

\newcommand{\data}{data/}

\title{Optimization under rare chance constraints\thanks{Resubmitted Jan 2022.
\funding{This material was based upon work
supported by the U.S. Department of Energy, Office of Science,
Office of Advanced Scientific Computing Research (ASCR) under
Contract DE-AC02-06CH11347. We acknowledge partial support by the US
National Science Foundation (NSF) through grants DMS \#1723211 to Shanyin Tong.}}}

\author{Shanyin Tong\thanks{Courant Institute, New York University, New York, NY (\email{shanyin.tong@nyu.edu}).}
\and Anirudh Subramanyam\thanks{Argonne National Laboratory, Lemont, IL (\email{asubramanyam@anl.gov}, \email{vhebbur@anl.gov}).}
\and Vishwas Rao\footnotemark[3]}

%% file: introduction.tex
\section{Introduction}

Rare high-impact events are phenomena that occur with low probability, but can have widespread--and sometimes even catastrophic--impacts on the economy, society and environment.
Examples of such phenomena include cascading blackouts in the power grid, stress fractures and material failures in bridges and other civil structures, extreme oceanic turbulence leading to natural disasters like hurricanes and tsunamis, financial stock market crashes, and large-scale pandemics like the ongoing COVID-19 outbreak.
It is of paramount importance to \emph{quantify the risk} of such rare high-impact events in a given system, by estimating their probability of occurrence, and to \emph{actively reduce this risk} by tuning the controllable properties of the system. %

Chance constraints provide an attractive modeling framework to mitigate the risk of such rare high-impact phenomena, since it is often difficult--or even impossible--to estimate the costs associated with the underlying events.
This motivates our study of chance constrained optimization problems of the following form.
\begin{equation}\label{eq:ccp}
\begin{aligned}
& \underset{u \in \mathcal{U}}{\text{minimize}}
& & J(u) \\
& \text{subject to}
& & \PP\left(F(u, \xi) \geq z\right)\leq \alpha,\quad \text{where } \alpha \ll 1.
\end{aligned}
\end{equation}
Here, $u$ represents the set of control or decision variables with domain $\mathcal{U}$, $J : \mathcal{U} \to \mathbb{R}$ represents the cost associated with decisions $u$, and $\xi$ denotes uncertain parameters that can influence the decisions via the function $F$.
Thus, $\mathcal{U}$ represents deterministic constraints on our decisions that are unaffected by uncertainties, and we consider $\mathcal{U} \subseteq \mathbb{R}^m$ or $\mathcal{U} \subseteq L^2(\mathcal{D})$ for some $\mathcal{D}\subseteq \mathbb{R}^m$.
We shall assume that the uncertain parameters $\xi$ are modeled as random vectors with probability distribution $\Pm$ supported on $\Xi \subseteq \mathbb{R}^n$,
and use $\PP(A)$ to denote the probability of an event $A$ under the measure induced by $\Pm$.
Therefore, $F : \mathcal{U} \times \Xi \to \mathbb{R}$ represents the uncertainty-affected system metric (also variously referred to as the \emph{parameter-to-observation map} or \emph{limit-state function}),
and $z$ is a specified upper bound for the metric $F$ such that, for any fixed decision vector $u$ and fixed realization $\xi$ of the random vector, $F(u, \xi) \geq z$ leads to an ``unsafe condition'' or ``bad event'' that we seek to mitigate.

Problem~\eqref{eq:ccp} thus seeks to determine decisions $u$ such that the probability of the system being unsafe, $\PP\left(F(u, \xi) \geq z\right)$, remains below a certain specified risk threshold $\alpha \in (0, 1)$.
We are particularly interested in the \textbf{rare event regime} where, for reasonable values of $u \in \mathcal{U}$, the latter quantities--the probability of being unsafe and the risk threshold $\alpha$--are very close to $0$.
We defer further assumptions and concrete examples to the subsequent sections.

One of the main challenges in solving \eqref{eq:ccp} is to tackle the chance constraint for extreme values of the risk threshold $\alpha\ll 1$ or the upper bound $z \to \infty$, since the resulting probability $\PP\left(F(u, \xi) \geq z\right)$ can be very small in these cases.
For general nonlinear $F$, classical methods %
approximate this probability by generating samples from the original distribution of (or some other `trial' distribution closely related to) $\xi$. %
However, the rare nature of the event coupled with the potential nonlinearity and complexity of $F$ render sampling-based methods generally intractable and potentially inaccurate. 
As \Cref{fig:sample} illustrates, such samples concentrate around regions of high likelihood, and relatively few samples cover the unsafe region.
Indeed, assuming that $F$ satisfies certain convexity properties, the number of samples required to obtain a feasible solution to~\eqref{eq:ccp} is $\mathcal{O}(\alpha^{-1})$ (\textit{e.g.}, see~\cite{calafiore2005uncertain}), which becomes impractical when $\alpha$ is $10^{-6}$ or $10^{-10}$ (say).
When taking into account problem-dependent constants as well, the actual required sample size (even for modest $\alpha$ and Gaussian $\xi$) can become arbitrarily large~\cite[Example~1]{Henrion2013:sample_size}.
This is made doubly worse when these samples are included in the optimization problem, since $F$ must be evaluated at each generated sample and this can lead to computational difficulties if $F$ itself is complex.
This motivates the development of an alternate \emph{sampling-free} method to solve \eqref{eq:ccp}, whose computational complexity does not scale with extreme values of $z$ or $\alpha$.

\input{Samples_Illustration.tex}

\subsection{Overview of methodology and contributions} 

We propose a sampling-free method to solve~\eqref{eq:ccp} by exploiting--and extending--recently developed results based on large deviation theory (LDT).
Under certain assumptions on the distribution of $\xi$ and the function $F$, it was shown in~\cite{dematteis2019extreme,tong2020extreme} that, for fixed $u$ and $z$, the rare event probability $\PP(F(u,\xi) \geq z)$ can be estimated using only local derivative information of $F$ evaluated at the \emph{LDT minimizer}:
\begin{equation}
\label{eq:LDT-min}
\xistar \in \argmin_{\xi \in \Xi} \left\{I(\xi) : F(u, \xi) \geq z\right\},
\end{equation}
where $I$ is the so-called \emph{rate function} that depends only on the distribution of $\xi$.
Although the method is more generally applicable, we focus our attention on the case where $\xi$ follows a general Gaussian mixture distribution.
Depending on the extent of local derivative information of $F$ that one can exploit or access, we then propose first- and second-order approximations of the rare event probability by quantifying the measure of the corresponding set bounded by the first- and second-order Taylor expansions of $F$ at $\xistar$,
and provide easily verifiable conditions under which the proposed estimates are exact or constitute upper bounds on the true probability.
However, using these estimates in place of the true probability in \eqref{eq:ccp} amounts to a bilevel optimization problem, since the LDT minimizer \eqref{eq:LDT-min} is a function of the decision variables $u$; to tackle this, we replace \eqref{eq:LDT-min} with its optimality conditions, and elucidate when the latter are sufficient to characterize the minimizer.
To demonstrate the effectiveness of our method, and its potential benefits over classical sampling-based methods, we consider a variety of applications from different domains. %
Specifically, via problems in financial portfolio management, structural column design and optimal control of a system of partial differential equations (PDE), we illustrate the accuracy of our probability estimates, feasibility of the obtained optimal solutions, improvement in computation time, and ability to handle complex nonlinear models.

The main contributions of our work may be summarized as follows.
\begin{enumerate}
    \item We provide explicit formulas for the first- and second-order LDT approximations of rare event probabilities when the constraint function $F$ is a continuously differentiable function of random parameters that follow a general Gaussian mixture distribution. This generalizes existing results~\cite{tong2020extreme} for Gaussians, and can be used to tackle general continuous distributions by approximating them with a Gaussian mixture. We also provide an alternate characterization of the second-order LDT approximation for Gaussians that is more numerically efficient for our purposes. %
    \item We apply the LDT approximations to optimization problems affected by rare chance constraints, and provide reformulations that can be solved by off-the-shelf solvers. Notably, our formulations are sampling-free and scalable: their size and complexity does not scale with the rarity of the chance constraints. We also contribute easily verifiable conditions under which their optimal values are exact or provide conservative upper bounds to the true optimal values. Finally, we do not limit ourselves to linear or convex $F$ and only require that it is twice differentiable, %
    thus enabling fast local solutions to nonlinear and nonconvex problems.
    \item We conduct computational experiments on linear, nonlinear and PDE-constrained optimization problems from several application areas to demonstrate the viability of our method and compare its performance with traditional sampling-based methods.
\end{enumerate}

This paper is structured as follows.
After introducing some notation and basic assumptions in \cref{sec:notation}, we provide a literature review of solution methods for chance constrained problems in \Cref{sec:sample}.
In \Cref{sec:LDT}, we provide an overview of our techniques, the large deviation results we use, as well as easily verifiable sufficient conditions under which we can establish optimality and feasibility of our formulations; notably, this section does not make specific assumptions on the nature of the distributions.
In \Cref{sec:G} and \Cref{sec:GM}, we then specalize our results to Gaussian and Gaussian mixture distributions, respectively, and provide closed-form expressions of the probability estimates and optimization formulations.
In \Cref{sec:app}, we illustrate the computational performance of our method on several applications, and finally, in \Cref{sec:conclusions}, we provide concluding remarks and directions for future work.

\subsection{Notation and setup}\label{sec:notation}
We use $\RR_{+}$ to denote the set of non-negative reals, and $\mathrm{I}_n$ to denote the $n\times n$ identity matrix.
We define $[x]_+ \coloneqq \max\{x,0\}$, and $\indicator{X}$ to be the indicator function of set $X$: $\indicator{X}(x)$ is equal to $1$ if $x \in X$ and $0$ otherwise.
Given a symmetric matrix $Q \in \RR^{n\times n}$,
we use $Q \succ 0$ ($\succeq 0$) to denote that it is positive definite (positive semidefinite),
and in such cases, we denote the corresponding weighted inner product between $x \in \RR^n$ and $y \in \RR^n$ by
$\langle x, y \rangle_Q \coloneqq \langle x, Qy \rangle$,
and the induced norm by $\|x\|_Q$.
We denote the standard Euclidean inner product interchangeably using both $\langle x, y \rangle$ and $x^\top y$.

With a slight abuse of notation, we use the same symbol $\xi$ to denote the random parameter vector as well as its realizations; the meaning should be clear from the context.
We use $\Pm$  to denote the distribution of $\xi$, $\Xi$ to denote its support, and $\PP(A)$ to denote the probability of an event $A$ under the measure $\Pm$.
Given a measurable function $f(\xi)$, %
we let
$\EE(f(\xi))$ denote the expectation of $f(\xi)$.
Finally, we use $S: \RR^n \to \RR$ to denote the cumulant generating function:
\begin{equation}
\label{eq:cgf}
S(\eta) \coloneqq \log \EE e^{\<\eta,\xi \>} = \log \left\lbrace \int_{\Xi} e^{\<\eta,\xi \>} \mathrm{d}\Pm(\xi)\right\rbrace,
\end{equation}
and define the rate function $I: \Xi \to \RR$ as the convex (Fenchel) conjugate of $S$:
\begin{equation}
\label{eq:I}
I(\xi) \coloneqq \max_{\eta\in\RR^n} \left\{ \<\eta,\xi \> - S(\eta) \right\}.
\end{equation} 

Throughout the paper, we assume that 
$\mu \coloneqq \EE(\xi)$ exists and that $S$ is continuously differentiable on $\RR^n$.
These are easily verified to be true whenever $\Pm$ is a Gaussian (or a mixture of Gaussians). %
In particular, when working with multivariate Gaussian parameters in $\RR^n$, %
we use $\xi\sim\Gauss{\mu}{\Sigma}$ to denote that
$\xi$ is a multivariate Gaussian with mean $\mu\in
\RR^n$ and covariance matrix $\Sigma \succ 0$, and
$\xi \sim \sum_{i = 1}^M w_i \Gauss{\mu_i}{\Sigma_i}$ to indicate that $\xi$ follows a Gaussian mixture distribution with $M$ components, where component~$i \in \{1, \ldots, M\}$ is a Gaussian with mean $\mu_i$, covariance $\Sigma_i \succ 0$ and weight $w_i > 0$ such that $\sum_{i=1}^M w_i = 1$.
Finally, we use $\Phi$ to denote the cumulative distribution function (CDF) of the standard normal distribution $\mathcal{N}(0, 1)$. %

Regarding the structure of $F$, we shall assume that
$F(\cdot, \xi)$ is twice continuously differentiable on $\mathcal{U}$ for every $\xi \in \Xi$
and
$F(u, \cdot)$ is twice continuously differentiable on closure of $\Xi$ for every $u \in \mathcal{U}$.
In addition, we assume that
there exist $\emptyset \neq \mathcal{U}_0 \subset \mathcal{U}$, $z_0 < \infty$ and $K_0 > 0$ such that,
for all $u \in \mathcal{U}_0$ and all $z\in [z_0, \infty)$,
we have
$F(u, \mu) < z_0$,
$\Xi(u, z) \coloneqq \left\{ \xi \in \Xi : F(u, \xi) = z \right\} \neq \emptyset$ and
$\| \nabla_{\xi} F(u, \xi) \| \geq K_0 > 0$ for all $\xi \in \Xi(u, z)$.

We note that the existence of $u \in \mathcal{U}$ and $x$ in a neighborhood of $z$ for which $\Xi(u, x) \neq \emptyset$ is not restrictive;
indeed, if one has $F(u, \xi) < z$ for all $\xi \in \Xi$ and $u \in \mathcal{U}$, then \eqref{eq:ccp} is trivially feasible for all $u \in \mathcal{U}$ (\textit{i.e.}, the chance constraint can be ignored), and similarly it is trivially infeasible if $F(u, \xi) > z$ for all $\xi \in \Xi$ and all $u \in \mathcal{U}$.
The assumption $\| \nabla_{\xi} F(u, \xi) \| > 0$ ensures that the hypersurface $\Xi(u, x) \subset \Xi$ does not degenerate and remains smooth and simply connected,
and that the minimization problem defined in~\eqref{eq:LDT} is always feasible for any $u \in \mathcal{U}_0$. %
We defer further assumptions to the sections where they are used.

Finally, note that we only consider the case of a single \emph{individual chance constraint}. 
The extension to multiple individual chance constraints is straightforward and we do not present it for ease of exposition.
On the other hand, the extension to \emph{joint chance constraints} where $F: \mathcal{U} \times \Xi \to \mathbb{R}^p$ and $z \in \mathbb{R}^p$, with $p > 1$ is not straightforward.
In such cases, one can either resort to a Bonferroni approximation (\emph{i.e.}, union bound) to convert the joint chance constraint into $p$ individual chance constraints, or apply our method to a smooth approximation of the function
$
\phi(x, u) = \min \{F_1(u, \xi) - z_1, \ldots, F_p(u, \xi) - z_p\}
$
(\textit{e.g.}, see \cite{pena2020solving}).

%% file: Samples_Illustration.tex
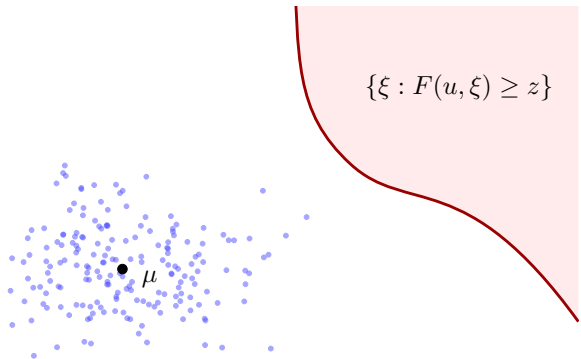
\begin{figure}[tbhp]\centering
    \begin{tikzpicture}[scale =0.9]
    \begin{axis}[compat=1.11, width=10cm, height=7cm,
    xmin=-6,
    xmax=4,
    ymin=-4,
    ymax=3,
    axis line style={draw=none},
    tick style={draw=none},
    yticklabels={,,},
    xticklabels={,,},
    ]
    \addplot [color=blue!70, only marks, mark=*,mark size=1pt,
    opacity=0.5]
    table[x=x,y=y] {\data/B0sample.txt};
    \filldraw [fill=red!8,draw=none]
    (4,-3) .. controls (2,0) and (1,-1)
    .. (0,0).. controls (-1,1) and (-0.9,2) .. (-1,4) -- (4,4)  -- (4,-3);
    \draw[red!60!black,very thick]  (4,-3) .. controls (2,0) and (1,-1)
    .. (0,0).. controls (-1,1) and (-0.9,2) .. (-1,4) ;
    \draw[black, yshift=4.5cm,xshift=8cm] node {$\{\xi: F(u,\xi)\geq z\}$};
    \addplot[mark=*] coordinates {(-4,-2)}  node[fill=white,yshift=-0.15cm,xshift=0.4cm]{$\mu$};
    \end{axis}
    \end{tikzpicture}
    \caption{Drawbacks of sampling-based methods for estimating rare events. The blue dots represent samples from the original distribution $\Pm$ with mean $\mu$, whereas the shaded region represents the rare event set $\{\xi: F(u,\xi)\geq z\}$ for some fixed decision vector $u$. Observe how very few samples fall inside the rare event set. %
    }\label{fig:sample}
\end{figure}

%% file: literature.tex
\section{Literature review} \label{sec:sample}

Chance constrained problems have been intensely studied for over five decades~\cite{shapiro2014lectures,prekopa2003}.
They are generically NP-hard since the set of feasible decisions is usually nonconvex, and admit 
tractable reformulations only in 
restrictive settings (\textit{e.g.}, 
when 
$F$ is bi-affine and $\xi$ follows a normal distribution~\cite{prekopa2003}).
Therefore, the majority of existing research takes one of two approaches: 
(1) exploit or assume some special structure in $F$ or the distribution of $\xi$, including linear or convex $F$ \cite{Lagoa2005,nemirovski2007convex}, finitely supported distributions of $\xi$ \cite{luedtke2010integer},
\edits{and elliptically symmetric Gaussian-like distributions of $\xi$ \cite{vanAckooij2014},}
or (2) resort to 
some form of 
sampling to approximate the chance constraint;
we discuss these approaches in \cref{sec:lit_review_sampling_methods}.
In \cref{sec:lit_review_reliability}, we discuss connections of our 
method to 
so-called 
reliability-based methods
that are prevalent 
in engineering design.

\subsection{Sampling-based methods}\label{sec:lit_review_sampling_methods}

Importance sampling is 
the most commonly used sampling method to estimate rare event probabilities~\cite{Kahn_1953}.
These methods belong to the class of variance reduction techniques, and sample from a (typically problem-specific) biasing distribution that is constructed so that probability mass is concentrated in the region of interest (\textit{e.g.}, see~\cite{Dunn_2011B, Asmussen_2007B, Bucklew_2013B}).
\edits{Most of these methods construct the biasing distribution for a fixed decision vector $u$.
	Therefore, they are not directly applicable to chance constrained optimization, where $u$ may change across algorithmic iterations.
	The quality of a single biasing distribution may not be uniform over these iterations, and this in turn may warrant the construction of multiple such distributions, albeit at the expense of additional computational effort.}
Nevertheless, existing methods for rare chance constrained optimization %
are also based on variants of importance sampling~\cite{rubinstein1997optimization,barrera2016chance
}.
\edits{For instance, \cite{barrera2016chance} constructs an importance sampling distribution that is uniformly efficient across the entire decision space; however, this is achieved by making the (restrictive) assumption that the components of $\xi$ are mutually independent random variables.}
Given $N$ independent samples or realizations, $\xi^1, \ldots, \xi^N$, that have been obtained via importance sampling (or some other method),
the probability, $\PP\left(F(u, \xi) \geq z\right)$ in \eqref{eq:ccp}, can then be approximated via one of two broad classes of approaches which we discuss in %
\cref{sec:lit_review_indicator} and \cref{sec:lit_review_scenario}.
Examples of generic sampling-based solution methods that do not fall into these classes can be found in~\cite{curtis2018sequential,pena2020solving}.
\edits{An alternative to importance sampling is quasi Monte Carlo sampling of the unit hypersphere that arises in the spherical-radial decomposition method for elliptically symmetric distributions, which we discuss in \cref{sec:lit_review_srd}.}
However, 
none of these approaches can be expected to work well in general for rare chance constraints, since the required number of samples can be $O(\alpha^{-1})$ in the worst case. 
Indeed, even if we manage to reduce the required number of samples \edits{using importance sampling or by spherical radial decomposition},
the former still entail modeling additional (at least $\mathcal{O}(N)$) variables and constraints, \edits{whereas the latter entail high per-iteration cost of constraint and gradient evaluations} in the optimization problem, which can make its solution computationally difficult.

\subsubsection{Indicator formulations}\label{sec:lit_review_indicator}
The methods in this class are all based on reformulating the probability in \eqref{eq:ccp} as the expectation of an indicator function
and then replacing the expectation with a \emph{sample average approximation} (SAA)~\cite{luedtke2008sample,pagnoncelli2009sample}:
\begin{equation*}
\label{eq:PMC}
\PP\left(F(u, \xi) \geq z\right)
=
\EE\left[\indicator{[z,\infty)}\left(F(u,\xi)\right)\right]
\approx
\frac{1}{N} \sum_{i=1}^{N} \indicator{[z,\infty)}\left(F(u,\xi^i)\right).
\end{equation*}
The discontinuous indicator function $\indicator{[z,\infty)}$ is then either exactly reformulated using $\mathcal{O}(N)$ binary variables (e.g., see~\cite{luedtke2010integer}) resulting in a mixed-integer nonlinear optimization problem,
or approximated using continuous functions to avoid introducing binary variables.
Examples of the latter include difference-of-convex approximations~\cite{hong2011sequential}, sigmoidal or other smooth approximations~\cite{geletu2017inner,adam2018solving,cao2020sigmoidal}
as well as
convex outer approximations~\cite{nemirovski2007convex}.
An example of a smooth approximation (that we also compare in our numerical experiments)
is the sigmoidal approximation proposed in~\cite{cao2020sigmoidal}.
It is based on the observation that, for any $\nu, \tau>0$, one can construct an outer-approximation of the indicator function that is asymptotically optimal as $\nu,\tau\to \infty$.
\begin{equation*}
\indicator{[z,\infty)}(x) \leq\left[\frac{2(\nu+1)}{\nu+\exp(-\tau(x-z))} -1\right]_+,
\quad
\lim_{\nu,\tau\to \infty} \left[ \frac{2(\nu+1)}{\nu+\exp(-\tau(x-z))} -1\right]_+ = \indicator{[z,\infty)}(x).
\end{equation*}
This allows us to obtain the following approximation to~\eqref{eq:ccp}:
\begin{equation}
\label{eq:op-SA2}
\begin{aligned}
& \underset{u \in \mathcal{U}, p \in \mathbb{R}^N_{+}}{\text{minimize}}
& & J(u) \\
& \text{subject to}
& & p_i\geq \frac{2(\nu+1)}{\nu+\exp(-\tau(F(u,\xi^i)-z))} -1 , \quad i \in \{1,\ldots,N \}, \\
& & & \frac{1}{N}\sum_{i=1}^{N} p_i \leq \alpha.
\end{aligned} \quad \text{(SAA)}
\end{equation}

Another class of general outer approximations was proposed in~\cite{nemirovski2007convex} by using any non-decreasing convex function $\psi : \RR \to \RR_{+}$ satisfying $\psi(0) = 1$ and $\psi(x) > 1$ for all $x > 0$.
In particular, by noting that
$
\indicator{[z,\infty)}\left(x\right) \leq \psi\left( t^{-1}\left( x - z \right) \right)
$
is satisfied for any such function and scalar $t > 0$,
it can be shown that satisfaction of the inequality,
\begin{equation*}
\inf_{t > 0}\left\lbrace t \EE\left( \psi\left( t^{-1}\left( x - z \right) \right) \right) - t\alpha \right\rbrace  \leq 0,
\end{equation*} 
implies satisfaction of the probability constraint in \cref{eq:ccp}.
Under the choice $\psi(x) = [1 + x]_{+}$, the above inequality is equivalent to
$
\CVaR_{1-\alpha}(F(u,\xi)-z)\leq 0,
$
where, for any measurable function $f(\xi)$ and any $\alpha \in (0, 1)$, $\CVaR\nolimits_{1-\alpha}(f(\xi)) \coloneqq \inf_{t\in\RR} \left\lbrace t + \alpha^{-1}\EE\left[f(\xi) - t \right]_+\right\rbrace$ is the conditional-value-at-risk (CVaR) of the random variable $f(\xi)$.
By approximating the expectation in the expression for CVaR with a sample average, we obtain the following conservative approximation to \eqref{eq:ccp}:
\begin{equation}
\label{eq:op-CVaR}
\begin{aligned}
& \underset{u \in \mathcal{U}, p \in \RR_{+}^N, t \in \RR}{\text{minimize}}
& & J(u) \\
& \text{subject to}
& & p_i\geq F(u,\xi^i)-z + t, \quad i \in \{1,\ldots,N \}, \\
& & & \frac{1}{N}\sum_{k=1}^{N} p_i \leq \alpha t.
\end{aligned}\quad \text{(CVaR)}
\end{equation}
The CVaR approximation has the advantage of preserving any convexity structure in $F$ (with respect to $u$), and so we use it as a benchmark in our numerical experiments.

\subsubsection{Scenario approaches}\label{sec:lit_review_scenario}
An alternative approach is to enforce the constraints deterministically for each of the $N$ realizations of $\xi$:
$
F(u,\xi^i)\leq z,
$
for all $i \in \{1, \ldots, N\}$.
As in the case of the CVaR approximation, this preserves any convexity structure in $F$ and the hope is that for sufficiently large $N$, optimal solutions of the scenario approximation will be feasible in \eqref{eq:ccp} with high probability.
Under the assumption that $F(\cdot, \xi)$ is convex for all $\xi \in \Xi$, the work~\cite{calafiore2005uncertain} shows that choosing $N = \mathcal{O}(m\alpha^{-1} \log(\alpha^{-1}))$ suffices, where $m$ is the number of decision variables $u$.
This size can be prohibitively large for small values of $\alpha$.
Therefore, \cite{nemirovski2006scenario,blanchet2020optimal} consider extensions of this approach to reduce the required sample size by making further assumptions on $F$ and the distribution of $\xi$.

\subsubsection{Spherical-radial decomposition}\label{sec:lit_review_srd}
\edits{
	Generic Monte Carlo sampling over the entire sample space can be inefficient when dealing with rare chance constraints, since the majority of samples will lie outside the rare event set.
	The spherical-radial decomposition (SRD) is an alternative approach that can achieve significant variance reduction whenever $\xi$ follows an elliptically symmetric distribution. %
	For example, when $\xi \sim \Gauss{\mu}{\Sigma}$ follows a Gaussian distribution with a Cholesky factorization $\Sigma=LL^\top$, one can write $\xi=\mu+r L \zeta$, where $\zeta$ is a random direction with uniform distribution $\nu$ on the unit sphere  $\mathbb{S}^{n-1} \coloneqq \{x \in \RR^n: \|x\| = 1\}$, and $r$ is the random magnitude of that direction following a chi-distribution with degree $n$.
	This allows us to express $\PP(F(u,\xi)\leq z)=\int_{\zeta\in \mathbb{S}^{n-1}} m_{\chi, n}(R(u,\zeta)) \mathrm{d} \nu \approx\frac{1}{N}\sum_{i=1}^N m_{\chi, n}(R(u,\zeta^i))$, where $m_{\chi, n}$ is the probability measure of the chi-distribution with degree $n$, $R(u,\zeta) = \{r>0: F(u, \mu + rL\zeta)\leq z \} $ are the magnitudes of direction $\zeta$ which intersect the complement of the rare event set $\{\xi: F(u,\xi)\leq z\}$, and $\zeta^1, \ldots, \zeta^N$ are uniform samples on the unit sphere $\mathbb{S}^{n-1}$ (e.g., obtained via quasi Monte Carlo sampling).
	Notably, \cite{vanAckooij2014,van2018sub} have shown that the same set of samples can also be used to compute the gradient (with respect to $u$) of $\PP(F(u,\xi)\leq z)$.
	\modifications{The crucial step in SRD is multiple root-finding to determine $R(u,\zeta)$, which 
		can become computationally expensive whenever explicit closed-form solutions of the roots (of $F(u, \mu + rL\zeta) = z$ in $r$ for fixed $u$ and $\zeta$) may not be available.
		This is especially true in problems where $F$ depends nonlinearly on $\xi$ (for e.g., via a PDE) and a good initial guess for $r$ may not be readily available.
		Moreover, this has to be repeated for each sampled direction $\zeta$ and each candidate decision $u$ that will change over optimization iterations. %
		Nevertheless, the method has found applications in several problems,
		including those involving joint chance constraints, where $F$ may be a vector-valued function.
		Successful applications include the probabilistic capacity maximization of gas networks \cite{heitsch2019probabilistic} %
		and optimal Neumann boundary control of a vibrating string with uncertain initial data and probabilistic terminal constraints \cite{farshbaf2020optimal}.
	}
}

\subsection{Reliability-based methods}\label{sec:lit_review_reliability}

In engineering, the structural reliability of buildings or bridges is assessed by solving a probability estimation problem similar to the one we consider, and reliability is enforced as part of the structural design by solving a so-called inverse reliability problem\cite{der1994inverse} or reliability-based design optimization problem \cite{valdebenito2010survey} that is similar in spirit to the chance constrained problem \eqref{eq:ccp}.
These methods approximate the distribution of $\xi$ using a (standard) multivariate normal distribution, and then compute the so-called most probable point (\textit{i.e.}, the point with maximum probability density) over the failure region $\{\xi \in \Xi: F(u, \xi) \geq z\}$ (\textit{e.g.}, see~\cite{du2001most}).
The distance of the most probable point from the origin is referred to as the reliability index,
which is not only used for estimating probabilities,
but also explicitly constrained as part of the optimal design problem.
In particular, the probability estimates are obtained by employing first- or second-order Taylor expansions of the limit-state function around the most probable point and then integrating the Gaussian densities over the set bounded by the Taylor approximations, resulting in the so-called first- and second-order reliability methods (FORM and SORM), respectively~\cite{schueller1987critical,du2001most}.
In contrast to a chance constraint on the limit-state function $F$, the reliability-based design problem explicitly constrains the reliability index and is solved using specialized procedures; e.g., see~\cite{lim2014second,valdebenito2010survey}. 

These methods share similarities with our LDT approach, especially for Gaussian distributions. In particular, both methods are sampling-free, and connect the probability estimations to an optimization problem. A crucial difference is that whereas reliability-based methods focus on finding the point with maximal likelihood, our LDT approach uses the point which characterizes the asymptotic behavior of the rare event probability, and thus includes information beyond just local probability densities. Another caveat is that reliability-based methods typically require transformation of the distribution to a standard normal, which can be difficult for a general high-dimensional distribution, whereas our methods are generally applicable for distributions with well-defined rate functions, including Gaussian mixtures that can, in limit, approximate arbitrary continuous distributions.
Finally, an important conceptual difference lies in the formulation of the optimization problem: whereas reliability-based methods %
use the reliability index as a surrogate for the probability to enforce constraints,
we treat the more traditional chance constrained formulation \eqref{eq:ccp}, and use the LDT minimizer along with local derivative information directly to approximate the probability in the chance constraint.

%% file: large_deviation.tex
\section{Large deviation theory for rare chance constraints} \label{sec:LDT}

Large deviation theory (LDT)~\cite{dembo1998large,varadhan1984large,touchette2011basic} is concerned with the asymptotic behavior of tails of probability distributions, and specifically the rates of exponential decay of probability measures of extreme events.
Classical LDT enables characterizing the asymptotic behavior of rare probabilities, such as the one in \eqref{eq:ccp}, with respect to some small (or large) parameter associated with the random variable $\xi$, while keeping $z$ fixed. It states that the behavior of the probability with respect to that small (or large) parameter is dominated by a decaying exponential term that can be characterized by the rate function of the random variable.
Recently, \cite{dematteis2019extreme} adapted classical theory and sharp asymptotics based on the dominating point of large deviations~\cite{iltis2000sharp,ney1983dominating,borovkov1965multi,broniatowski1995tauberian} to study the asymptotic behavior with respect to $z \to \infty$, while keeping the distribution of $\xi$ fixed. %
The key idea is to find a dominating point in the rare event set by solving an appropriately defined optimization problem and estimating probabilities solely using this dominating point.
In \cite{tong2020extreme}, it was shown that improved estimates can be obtained by exploiting Gaussianity as well as local approximations of $F$ around this dominating point.
The approach has been successfully applied to quantify rare probabilities of extreme oceanic phenomena~\cite{dematteis2018rogue,dematteis2019experimental,tong2020extreme}.

In the following subsections, we present an overview of this approach and show how it can be used for rare chance constrained optimization.
Specifically, in \cref{sec:prob}, we first reduce the problem of estimating the measure of the set bounded by (the potentially nonlinear) $F$ to estimating the measure of the set bounded by its first- and second-order Taylor expansions around a dominating point.
In \cref{sec:ldt_result}, we provide a specific choice of the dominating point based on large deviation arguments;
and finally, in Section \cref{sec:ldt_opt} we incorporate these probability estimates in the chance constrained problem.

\subsection{Taylor approximation of nonlinear limit state function}\label{sec:prob}
For a fixed choice of $u \in \mathcal{U}$ and $z \in \mathbb{R}$, a key challenge in computing the probability in \eqref{eq:ccp}, $\PP\left(F(u, \xi) \geq z\right)$ is the potential nonlinearity of $F$ with respect to $\xi$.
A simple way to tackle this is to construct first- and second-order Taylor approximations of $F$ around some point $\xistar \in \Xi$, %
\begin{subequations}\label{eq:Fapprox}
    \begin{alignat}{2}
    F_1(u,\xi; \xistar) \coloneqq& F(u,\xistar)+\<\nabla_\xi F(u,\xistar), \xi-\xistar \>, \label{eq:FFO}\\ 
    F_2(u,\xi; \xistar) \coloneqq& F(u,\xistar)+\<\nabla_\xi F(u,\xistar), \xi-\xistar \> +\frac{1}{2}\<\xi-\xistar, \nabla^2_\xi F(u,\xistar)( \xi-\xistar) \>. \label{eq:FSO}
    \end{alignat}
\end{subequations}
We can then compute first- and second-order probability estimates, $P_1$ and $P_2$ of $\PP\left(F(u, \xi) \geq z\right)$, by computing the measure of the sets bounded by the corresponding Taylor approximations:
\begin{equation}\label{eq:PFOPSO}
P_k(u,z,\xistar)
=
\Pm\left(\left\{\xi \in \Xi
: F_k(u, \xi; \xistar) \geq z \right\}\right), \quad k \in \{1, 2\}.
\end{equation}
It is important to note that this approach approximates the nonlinear (in $\xi$) chance constrained problem \eqref{eq:ccp} with a linear and quadratic (in $\xi$) chance constrained problem, respectively.
In particular, it allows us to exploit special structure (if any) in the distribution of $\xi$ to efficiently compute the first- and second-order approximations in \eqref{eq:PFOPSO}.
We provide examples of this computation for standard Gaussian and mixture distributions in the subsequent sections.
First, we make the immediate observation that $P_1$ provides a conservative approximation to the true probability whenever $F$ is concave in the uncertainty $\xi$.
\begin{proposition}[First-order probability estimate for concave functions]\label{prop:concave_LDT_prob_upper_bound}
    Fix $u \in \mathcal{U}$.
    Suppose that $\Xi$ is convex and the function $F(u, \cdot)$ is concave on $\Xi$.
    Then, $\PP\left(F(u, \xi) \geq z\right) \leq P_1(u,z,\xistar)$
    for any %
    $z \in \RR$ and $\xistar \in \Xi$.
\end{proposition}
\begin{proof}
    Since $F(u, \cdot)$ is continuously differentiable (refer to \cref{sec:notation}), it is concave on $\Xi$ if and only if
    \[
    F(u,\xi) \leq F(u,\xistar) + \<\nabla_\xi F(u,\xistar), \xi-\xistar\> = F_1 (u,\xi; \xistar), \quad \forall \xi \in \Xi.
    \]
    Therefore, we have the relationship
    \begin{equation*}
    \left\{\xi \in \Xi: F(u, \xi) \geq z \right\}\subseteq\left\{\xi \in \Xi: F_1 (u,\xi; \xistar) \geq z \right\},
    \end{equation*}
    \textit{i.e.}, the rare event set is included in the half-space bounded by $F_{1}(u,\xi; \xistar) = z$.
    Therefore,
    \begin{equation*}
    \PP\left(F(u, \xi) \geq z\right) =\Pm\left(\left\{\xi \in \Xi
    : F(u,\xi) \geq z \right\}\right)  \leq \Pm\left(\left\{\xi \in \Xi
    : F_{1}(u,\xi; \xistar) \geq z \right\}\right)=P_{1}(u,z,\xistar),
    \end{equation*}
    and the claim follows. %
\end{proof}

\input{FormSorm_Illustration.tex}

\subsection{Large deviation theory}\label{sec:ldt_result}
Although the Taylor approximations of $F$ provide a tractable means to estimate the true probability, their accuracy depends critically on the choice of $\xistar$.
In the following, we invoke large deviation results from \cite{dematteis2019extreme} to justify a choice of $\xistar$ that is asymptotically optimal as $z \to \infty$.
The formal statement of the result needs additional assumptions on the distribution of $\xi$ and the structure of $F$. The procedure for choosing $\xistar$ based on this result, and the first- and second-order approximations of $F$ at this point are shown in \Cref{fig:LDT}.
\begin{theorem}[Theorem~2.1~in~\cite{dematteis2019extreme}]\label{thm:ldt_main}
    Fix $u \in \mathcal{U}_0$. Suppose the following assumptions hold:
    \begin{enumerate}[label=(A\arabic*),leftmargin=*]
        \item
            \label{assume:unique_global_min}
            The optimization problem
            $
            \mathop{\min}_{\xi \in \Xi} \left\{I(\xi) : F(u, \xi) \geq z\right\}
            $
            has a unique global minimizer for all $z \in [z_0, \infty)$. %
            Therefore, %
            we can define:
            \begin{equation}
            \label{eq:LDT-min-opt}
            \xistar : \mathcal{U}_0 \times [z_0, \infty) \to \Xi, \quad \xistar(u, z) \coloneqq \argmin_{\xi \in \Xi} \left\{I(\xi) : F(u, \xi) \geq z\right\}.
            \end{equation}
        \item
            \label{assume:increasing_rate_function}
        	The function
        	$\xistar$ in \eqref{eq:LDT-min-opt} is continuously differentiable,
            and %
            there exists $K_0' > 0$ such that the rate function $I(\xistar(u, \cdot))$ in \eqref{eq:I} is strictly increasing with
        	\begin{equation}
        	\label{eq:LDT-I_increasing}
        	I(\xistar(u,z))\to \infty \quad \text{and} \quad \|\nabla
        	I(\xistar(u,z))\|\ge K_0' >0
        	\quad \text{as}\quad z\to\infty.
        	\end{equation}
        \item
            \label{assume:halfspace}
            The following relation holds for all $z \in [z_0, \infty)$:
            \begin{equation*}
            \left\{\xi \in \Xi: F(u, \xi) \geq z \right\}\subseteq\left\{\xi \in \Xi: \<\nabla_\xi F(u,\xistar(u, z)), \xi-\xistar(u, z)\> \geq 0 \right\},
            \end{equation*}
        \item
            \label{assume:finite_correction}
           	There exists $K_1 > 0$ such that, for all $u \in \mathcal{U}_0$ we have
           	\begin{equation*}
           	\label{eq:LDT-ass4}
           	\lim\limits_{z\to\infty}
               \frac{
                       \log\left( \int _{\Lambda(u, z, K_1) } \exp\left(
                       I(\xistar(u, z)) +
                   			\<\nabla I(\xistar(u, z)), \xi-\xistar(u, z)\> \right)
                   		\mathrm{d}\Pm(\xi)\right)
               }{
               I(\xistar(u, z))
           }
                = 0
           	\end{equation*}
               where
               $\displaystyle
               \Lambda(u, z, s) \coloneqq \left\{\xi \in \Xi :
               F(u, \xi) \geq z,
               \frac{\left<
                    \nabla_\xi F(u,\xistar(u, z)), \xi-\xistar(u, z)
               \right>}{\left\lVert \nabla_\xi F(u,\xistar(u, z)) \right\rVert}
                < s
                \right\}
               $.
    \end{enumerate}
    Then, the following result holds.
    \begin{equation}\label{eq:ldt_main}
    \lim\limits_{z\to\infty}\frac{\log \PP\left(F(u, \xi) \geq z\right)}{I(\xistar(u, z))} = -1.
    \end{equation}
\end{theorem}

\Cref{thm:ldt_main} implies that under its stated assumptions, the rare event probability in~\cref{eq:ccp} satisfies the asymptotic equality
\begin{equation}
\PP\left(F(u, \xi) \geq z\right) \asymp \exp\left(-I(\xistar(u, z))\right).  \label{eq:LDT}
\end{equation}
where $\asymp$ indicates that the ratio of logarithms of both sides tends to 1 as $z\to\infty$,
For large but finite values of $z$, the right-hand of~\cref{eq:LDT} can therefore be expected to serve as a good approximation of the true probability.
This justifies the LDT minimizer $\xistar(u, z)$ defined in \eqref{eq:LDT-min-opt} as the ideal point around which one should Taylor approximate the nonlinear $F$.

We now review the assumptions leading to the result of \Cref{thm:ldt_main} in detail, with the goal of simplifying them. %
Although not explicitly stated in \cite{dematteis2019extreme}, assumption~\ref{assume:unique_global_min} is necessary to ensure that the LDT minimizer in \cref{eq:LDT-min-opt} is well defined.
Assumption~\ref{assume:increasing_rate_function} formally states that the event $\{F(u, \xi) \geq z\}$ becomes rare as $z \to \infty$,
whereas assumption~\ref{assume:halfspace} requires that $\xistar$ is a so-called dominating point~\cite{ney1983dominating}: that is, the rare event set $\{\xi \in \Xi: F(u, \xi) \geq z\}$ is completely enclosed in the halfspace that is tangent to the hypersurface $\Xi(u, z)$ at $\xistar$.
Finally, Assumption~\ref{assume:finite_correction} is an equivalent rephrasing of Assumption~5 in \cite{dematteis2019extreme}.
Specifically, it ensures that the true rare event probability does not decay much faster than $\exp(-I(\xistar(u, z)))$ as $z \to \infty$; \emph{i.e.}, that the right-hand side of \cref{eq:ldt_main} is a lower bound on its left-hand side expression.

The assumptions leading to the result of \Cref{thm:ldt_main} may not always be straightforward to verify.
With a goal of obtaining easily verifiable conditions, the following theorem provides alternative sufficient conditions that lead to the same large deviation result as \Cref{thm:ldt_main}.
\begin{theorem}[Sufficient conditions for \Cref{thm:ldt_main}]\label{thm:ldt_verifiable}
    Fix $u \in \mathcal{U}_0$. Suppose that the following condition is satisfied:
    \begin{enumerate}[label=(C\arabic*),leftmargin=*]
        \item\label{assume:concave_F}
        $\Xi = \RR^n$ and the function $F(u, \cdot)$ is concave on $\Xi$.
    \end{enumerate}
    Then, the following inequality is valid:
    \begin{equation}\label{eq:ldt_main_ub}
    \lim\limits_{z\to\infty}\frac{\log \PP\left(F(u, \xi) \geq z\right)}{I(\xistar(u, z))} \leq -1.
    \end{equation}
    Additionally, suppose that the following conditions also hold:
    \begin{enumerate}[label=(C\arabic*),leftmargin=*]\setcounter{enumi}{1}
        \item\label{assume:gradF_lipshitz}
        $\nabla_{\xi} F(u, \cdot)$ is Lipschitz continuous on $\Xi$.
        
        \item\label{assume:finte_volume}
        For any finite $r > 0$, the measure $\mathcal{P}$ satisfies
        \begin{equation*}
            \log\left( \int _{\mathcal{B}(r, \hat{\xi})} \exp\left(
            I(\hat{\xi}) + 
            \<\nabla I(\hat{\xi}), \xi-\hat{\xi}\> \right)
            \mathrm{d}\Pm(\xi)\right)
        = o(
            I(\hat{\xi})
            )
        \quad
        \text{ as } \| \hat{\xi} \| \to \infty,
        \end{equation*}
        where $\mathcal{B}(r, \hat{\xi}) \subset \Xi$ denotes a closed ball of radius $r$ that contains $\hat{\xi} \in \Xi$.
    \end{enumerate}
    Then, \Cref{eq:ldt_main} holds.
\end{theorem}
\begin{proof}
    To show the first part, we will show that condition~\ref{assume:concave_F} implies assumptions \ref{assume:unique_global_min}, \ref{assume:increasing_rate_function} and \ref{assume:halfspace} stated in \Cref{thm:ldt_main}.
    The proof for the validity of \cref{eq:ldt_main_ub} then follows from the proof of Theorem~2.1~in~\cite{dematteis2019extreme}.
    By a similar argument, we will show that the conditions taken together also imply assumption \ref{assume:finite_correction}.
    The proof that the left-hand side of \eqref{eq:ldt_main_ub} is $\geq$ to its right-hand side then follows from similar arguments as Theorem~2.1~in~\cite{dematteis2019extreme}.
    
    \paragraph{Verification of \ref{assume:unique_global_min}}
    Since the cumulant generating function $S$ is finite on $\RR^n$ (see \cref{sec:notation}),
    \cite[Theorem~2]{pistone1999} implies that it is strictly convex and twice continuously differentiable.
    Therefore, the rate function $I$, which by definition~\cref{eq:I} is the convex conjugate of $S$, satisfies the following properties:
    \textit{(i)} $I$ is strictly convex and twice continuously differentiable on $\RR^n$ \cite[Corollary~4.2.10]{HiriartUrruty1993}, and
    \textit{(ii)} $I$ is coercive \cite[Theorem~11.8]{Rockafellar}, \textit{i.e.}, $I(\xi) \to +\infty$ as $\| \xi \| \to \infty$.
    Together with the concavity of $F$ from condition \ref{assume:concave_F}, these properties imply that 
    $
    \mathop{\min}_{\xi \in \Xi} \left\{I(\xi) : F(u, \xi) \geq z\right\}
    $
    is a convex optimization problem (which is always feasible for $z \in [z_0, \infty)$, see \cref{sec:notation})
    with a finite optimal value
    that is attained by a unique global minimizer $\xistar(u, z)$ defined in \cref{eq:LDT-min-opt}.
    Notably, this implies assumption \ref{assume:unique_global_min}.
    
    \paragraph{Verification of \ref{assume:increasing_rate_function}}
    We show the implication of assumption~\ref{assume:increasing_rate_function} in two parts.
    \begin{itemize}[leftmargin=*]
        \item First, note that continuity of $F(u, \cdot)$ implies that
        $
        \sup_{\xi \in \Xi} \left\{ F(u, \xi) : \| \xi \| \leq r_0\right\} < \infty
        $
        for all finite $r_0 > 0$.
        This, in turn, implies that any point in the set
        $
        \left\{\xi \in \Xi : F(u, \xi) \geq z\right\}
        $
        must satisfy $\| \xi \| \to \infty$ as $z \to \infty$.
        In particular, we have  $\| \xistar(u, z) \| \to \infty$ and from the coercivity of $I$ established in \textit{(ii)} above, this implies that
        $I(\xistar(u,z))\to \infty$.
        
        \item Second, for any $z \in [z_0, \infty)$, there always exists a point $\xi \in \Xi = \RR^n$ satisfying $F(u, \xi) > z$, see \cref{sec:notation}.
        This implies that the convex problem
        $
        \mathop{\min}_{\xi \in \Xi} \left\{I(\xi) : F(u, \xi) \geq z\right\}
        $
        satisfies Slater's constraint qualification,
        and strong Lagrangian duality holds:
        \[
        I(\xistar(u,z)) = \max_{\lambda \in \RR_{+}} \left\{ \lambda z + \min_{\xi \in \Xi} \left\{ I(\xi) - \lambda F(u, \xi) \right\} \right\}
        \]
        For sufficiently large values of $z$, the maximizer in the right-hand side of the above expression is always attained at some $\lambda \geq \lambda_0 > 0$; indeed, if it is attained at $\lambda = 0$, then $I(\xistar(u,z)) =  \min_{\xi \in \Xi} I(\xi) < +\infty$ which would contradict $I(\xistar(u,z))\to \infty$.
        Therefore, as a function of $z$, the right-hand side of the above expression, and hence $I(\xistar(u,z))$, is convex and strictly increasing.
        Moreover,
        $\|\nabla I(\xistar(u,z))\| = \lambda \| \nabla_{\xi} F(u, \xi) \| \geq \lambda_0 K_0 > 0$,
        (Note that, by the same reasoning, we also have $F(u, \xistar(u,z)) = z$ which will be used in the next paragraph).
    \end{itemize}
    
    \paragraph{Verification of \ref{assume:halfspace}}
    This follows from the definition of concavity in condition~\ref{assume:concave_F}.

    \paragraph{Verification of \ref{assume:finite_correction}}
    Let $L > 0$ denote the Lipschitz constant in condition \ref{assume:gradF_lipshitz}.
    Consider now the closed ball $B(r, \xistar(u, z))$ of radius $r = K_0L^{-1} > 0$ centered at $\xistar(u, z) + r\hat{n}(u, z)$, where $K_0$ was defined in \cref{sec:notation}, and $\hat{n}(u, z)$ is the unit vector along $\nabla_\xi F(u,\xistar(u, z))$.
    We claim that $B(r, \xistar(u, z)) \subseteq \Lambda(u, z, 2r)$, where $\Lambda$ is defined in assumption \ref{assume:finite_correction}.
    If true, then observe that, the fact that $\| \xistar(u, z) \| \to \infty$ as $z \to \infty$ (established above), along with condition \ref{assume:finte_volume} evaluated at $\hat{\xi} = \xistar(u, z)$ implies \ref{assume:finite_correction} is satisfied with $K_1 = 2K_0L^{-1} > 0$.
    To show $B(r, \xistar(u, z)) \subseteq \Lambda(u, z, K_1)$,
    we proceed as follows,
    dropping the dependence on $u$ and $z$ for simplicity of exposition.
    \begin{itemize}
        \item First, we show $\xi \in B(r, \xistar)$ implies that $\< \hat{n}, \xi - \xistar\> \leq K_1$.
        \begin{align*}
        \xi \in B(r, \xistar) \iff  \| \xi  - \xistar - \hat{n} r \| \leq r %
        \iff \| \xi  - \xistar \| \leq 2r
        \implies
        \< \hat{n}, \xi - \xistar\> \leq 2r = K_1.
        \end{align*}
        
        \item Next, we show that $\xi \in B(r, \xistar)$ implies that $F(\xi) \geq z$.
        Lipschitz continuity of $\nabla_{\xi} F(u, \cdot)$ along with its twice differentiability implies that $\nabla_{\xi}^2 F(u, \xi) \succeq -L\mathrm{I}_n$ holds for all $\xi \in \Xi$. %
        Along with Taylor's theorem, this implies
        \begin{align*}
        F(u,\xi) &\geq F(u,\xistar) + \< \nabla_{\xi}F(u,\xistar), \xi - \xistar \> - \frac{L}{2} \| \xi - \xistar \|^2 \\
        &= z + \| \nabla_{\xi}F(u,\xistar) \| \left[ \< \hat{n}, \xi - \xistar \> - \frac{L}{2\|\nabla_{\xi}F(u,\xistar) \|} \| \xi - \xistar \|^2 \right] \\
        &\geq z + \| \nabla_{\xi}F(u,\xistar) \| \left[ \< \hat{n}, \xi - \xistar \> - \frac{1}{2r} \| \xi - \xistar \|^2 \right] \\
        &= z \quad (\text{since } \xi \in B(r, \xistar))
        \end{align*}
        where we used the fact that $F(u,\xistar) = z$, $\| \nabla_{\xi}F(u,\xistar) \| \geq K_0$ and the definition of $r = K_0 L^{-1}$.
    \end{itemize}
    Thus, we have also verified assumption \ref{assume:finite_correction}.
\end{proof}

Unlike \Cref{thm:ldt_main}, it is crucial to note that the conditions in \Cref{thm:ldt_verifiable} either depend only on the structure of $F$ or on the measure $\mathcal{P}$, and hence, are relatively easier to verify in general.
For example, \Cref{thm:gaussian_mixture_finite_volume} in \Cref{sec:GM} shows that when $\mathcal{P}$ is a Gaussian or mixture of Gaussian distributions, then condition \ref{assume:finte_volume} is automatically satisfied.

\subsection{Incorporating large deviation estimates in the chance constrained problem}\label{sec:ldt_opt}

We can now incorporate in the chance constrained problem \eqref{eq:ccp}, the probability estimates obtained using the first- and second-order Taylor expansions of $F$ that we introduced in \cref{sec:prob}, around the dominating point $\xistar$ from \cref{sec:ldt_result}.
Depending on whether we use the first- ($k = 1$) or second-order ($k=2$) approximation, we obtain the following conceptual problem:
\begin{equation}\label{eq:op-bi}
\begin{aligned}
& \underset{u \in \mathcal{U}}{\text{minimize}}
& & J(u) \\
& \text{subject to}
&& P_k(u,z,\xistar) \leq \alpha  \\ %
&&& \xistar \in \argmin_{\xi \in \Xi} \left\{I(\xi) : F(u, \xi) \geq z\right\}
\end{aligned}
\end{equation}
Observe that, since the dominating point depends on the choice of the decision variables $u$, the resulting problem is a generic nonlinear bilevel optimization problem.
Moreover, even though we use Taylor approximations of $F$ in constructing $P_k$, see \eqref{eq:PFOPSO}, the lower-level problem that characterizes the dominating point
$\xistar$
utilizes the original nonlinear function.
This is crucial for the validity of the asymptotic results in \Cref{thm:ldt_main} and \Cref{thm:ldt_verifiable}.

To tackle the bilevel structure, we replace the lower-level problem with its first-order optimality conditions.
Although it is straightforward to consider the extension where the support $\Xi = \{\xi \in \RR^n: G_j(\xi) \leq 0, \; j \in \{1, \ldots, J\}\}$ has an explicit algebraic description in terms of inequalities,
we assume that $\Xi = \RR^n$ for simplicity of exposition.
We have the following single-level approximation of the chance-constrained problem \eqref{eq:ccp}:
\begin{equation}\label{eq:op-LDT}
\begin{aligned}
& \underset{u, \xistar, \lambda}{\text{minimize}}
& & J(u) \\
& \text{subject to}
&& u \in \mathcal{U}, \; \xistar \in \Xi, \; \lambda \in \RR_{+} \\
&&& P_k(u,z,\xistar) \leq \alpha, \\
&&& F(u,\xistar)=z,\\
&&& \nabla I(\xistar)  = \lambda \nabla_\xi F(u,\xistar).
\end{aligned}
\end{equation}
Note that we have replaced $F(u,\xistar)\geq z$ with equality.
Doing so relies on an implicit (yet mild) assumption that $F(u, \bar{\xi}) < z$, where $\bar{\xi} \coloneqq \argmin_{\xi \in \Xi} I(\xi)$, and we provide a formal argument in \Cref{thm:LDT-ccop} below.
In any case, note that $\bar{\xi}$ is unique and always exists (see proof of \Cref{thm:ldt_verifiable}), and continuity of $F$ implies that this assumption is expected to hold for sufficiently large values of $z$.
It is important to note that this argument does not rely on any convexity assumptions on $F$.

Finally, we note that the rate function $I$ may not always be explicitly available depending on the distribution of $\xi$.
In such cases, we can instead rely on the fact that $I$ is the conjugate of the cumulant generating function $S$, for which we typically have efficient statistical codes that can compute its values and gradients for several distributions.
Indeed, from the definition of the rate function \eqref{eq:I} and convexity and differentiability of $S$ (see \cref{sec:notation}), we have
\begingroup
\allowdisplaybreaks
\begin{equation}\label{eq:xistar-eta}
\begin{aligned}
\displaystyle \xistar & \displaystyle \in \argmin_{\xi \in \Xi} \left\{I(\xi) : F(u, \xi) \geq z\right\} \\
\displaystyle\iff
\exists \etastar \in \RR^n:
(\xistar, \etastar) &\displaystyle \in \left\{
\begin{array}{c@{\;\;}l}
\displaystyle \argmin_{\xi \in \Xi, \eta \in \RR^n} & \displaystyle\<\eta,\xi \>-S(\eta)\\
\displaystyle \text{subject to} & \displaystyle F(u,\xi) \geq z,\;\; \xi=\nabla_\eta S(\eta)
\end{array}
\right\}.
\end{aligned}
\end{equation}
\endgroup
It is straightforward to verify that we then have $\nabla I(\xistar) = \etastar$.
We then obtain the following single-level approximation of the chance-constrained problem \eqref{eq:ccp}:
\begin{equation}\label{eq:op-LDT-etastar}
\begin{aligned}
& \underset{u, \xistar, \etastar, \lambda}{\text{minimize}}
& & J(u) \\
& \text{subject to}
&& u \in \mathcal{U}, \; \xistar \in \Xi, \; \etastar \in \RR^n, \; \lambda \in \RR_{+} \\
&&& P_k(u,z,\xistar) \leq \alpha, \\
&&& F(u,\xistar)=z,\\
&&& \etastar  = \lambda \nabla_\xi F(u,\xistar), \\
&&& \xistar =  \nabla S(\etastar).
\end{aligned}
\end{equation}
We establish some properties of the single-level approximation \eqref{eq:op-LDT} or equivalently \eqref{eq:op-LDT-etastar}.
\begin{theorem}[Properties of the LDT chance-constrained optimization]\label{thm:LDT-ccop}
    Suppose that $\Xi = \RR^n$ and $F(u, \bar{\xi}) < z$ for all feasible solutions $u \in \mathcal{U}$ of \eqref{eq:ccp} and \eqref{eq:op-bi}, where $\bar{\xi} \coloneqq \argmin_{\xi \in \Xi} I(\xi)$.
    Then, the following statements are true.
    \begin{enumerate}[leftmargin=*]
        \item\label{opt:bilevel_to_single_level}
        If $(u, \xistar)$ is feasible in the bilevel problem \eqref{eq:op-bi}, then there exists $\lambda \in \RR_{+}$ such that
        $(u, \xistar, \lambda)$ is feasible in the single-level problem \eqref{eq:op-LDT}.
        
        \item\label{opt:bilevel_to_single_level_convave_F}
        Under condition \ref{assume:concave_F}, the converse is also true: if $(u, \xistar, \lambda)$ is feasible in \eqref{eq:op-LDT}, then $(u, \xistar)$ is feasible in \eqref{eq:op-bi}. Therefore, \eqref{eq:op-LDT} is an exact single-level reformulation of \eqref{eq:op-bi}.
        
        \item\label{opt:firstorder}
        Under condition \ref{assume:concave_F} and $k = 1$, the first-order LDT problem \eqref{eq:op-LDT} is a conservative approximation of the chance constrained problem \eqref{eq:ccp} for any $z \in \RR$; that is, if $(u, \xistar, \lambda)$ is feasible in \eqref{eq:op-LDT}, then  $u$ is feasible in \eqref{eq:ccp}.

        \item\label{opt:exact}
        Under conditions \ref{assume:concave_F}, \ref{assume:gradF_lipshitz} and \ref{assume:finte_volume}, the feasible solutions of both the first- ($k = 1$) and second-order ($k = 2$) LDT problems \eqref{eq:op-LDT} are in one-one correspondence with those of the chance constrained problem \eqref{eq:ccp} and therefore, the former constitutes an asymptotic reformulation of the latter, as $z \to \infty$.
    \end{enumerate}
\end{theorem}
\begin{proof}
    We consider each of the claims separately.
    \paragraph{Proof of claim \ref{opt:bilevel_to_single_level}}
    From \cref{sec:notation}, we know that $\| \nabla_{\xi} F(u, \xi) \| > 0$ for all $u \in \mathcal{U}_0$.
    Therefore, all local minimizers of the lower-level problem in \eqref{eq:op-bi} satisfy the linear independence constraint qualification.
    Therefore, if $(u, \xistar)$ is feasible in \eqref{eq:op-bi}, then $\xistar$ must satisfy the first-order optimality conditions of the lower-level problem for some $\lambda \in \RR_{+}$.
    If $\lambda = 0$, then the stationarity condition $\nabla I(\xistar)  = \lambda \nabla_\xi F(u,\xistar)$ would imply that $\xistar = \bar{\xi}$; the primal feasibility of the lower-level problem then implies that $F(u, \xistar) = F(u, \bar{\xi}) \geq z$ which leads to a contradiction since $F(u, \bar{\xi}) < z$ by construction.
    Therefore, we must have $\lambda > 0$ and complementarity slackness then requires that $F(u, \xistar) = z$.
    
    \paragraph{Proof of claim \ref{opt:bilevel_to_single_level_convave_F}}
    Under condition \ref{assume:concave_F}, the lower-level problem in \eqref{eq:op-bi} satisfies Slater's constraint qualification, see also the proof of \Cref{thm:ldt_verifiable}.
    Therefore, its first-order optimality conditions are necessary and sufficient.
    
    \paragraph{Proof of claim \ref{opt:firstorder}}
    This is a direct consequence of \Cref{prop:concave_LDT_prob_upper_bound} and claim \ref{opt:bilevel_to_single_level_convave_F}.
    
    \paragraph{Proof of claim \ref{opt:exact}}
    Suppose that $(u, \xistar, \lambda)$ is feasible in \eqref{eq:op-LDT}.
    Then, under conditions \ref{assume:concave_F}, \ref{assume:gradF_lipshitz} and \ref{assume:finte_volume}, \Cref{thm:ldt_verifiable} ensures that this $\xistar$ is precisely the unique global minimizer appearing in \eqref{eq:LDT-min-opt}.
    Moreover, for this $u$ and $\xistar$, the measure of the rare event set $ \Pm\left( \left\{\xi \in \Xi: F(u, \xi) \geq z\right\}\right) \asymp \exp(-I(\xistar))$ as $z \to \infty$.
    An application of the same theorem to $F_k(u,\xi; \xistar)$ (see \eqref{eq:FFO} and \eqref{eq:FSO}) instead of $F$,
    must require that $\xistar$ again be the unique global minimizer appearing in \eqref{eq:LDT-min-opt};
    this is because of $\nabla_{\xi} F(u, \xistar) = \nabla_{\xi} F_k(u, \xistar; \xistar)$  (by construction) along with the strict convexity of the rate function $I$.
    Therefore, again we have $ \Pm\left( \left\{\xi \in \Xi: F_k(u, \xi; \xistar) \geq z\right\}\right) \asymp \exp(-I(\xistar))$ as $z \to \infty$.
    The claim now follows from claim \ref{opt:bilevel_to_single_level_convave_F}.
\end{proof}

We remark that the single-level reformulations \cref{eq:op-LDT} and \cref{eq:op-LDT-etastar} are \textit{conceptual} in the sense that they require an instantiation of the probability estimate $P_k$ as well as gradients of the rate function $I$ and/or the cumulant generating function $S$.
In \Cref{sec:G} and \Cref{sec:GM}, we provide explicit single-level formulations in the case where $\xi$ is a Gaussian and a mixture of Gaussians, respectively.

%% file: FormSorm_Illustration.tex
\begin{figure}[tbhp]\centering
	\begin{tikzpicture}[]
\begin{axis}[compat=1.11, width=9cm, height=7cm,
xmin=-5,
xmax=4,
ymin=-3,
ymax=4,
axis line style={draw=none},
tick style={draw=none},
yticklabels={,,},
xticklabels={,,},
]
\draw [blue,very thick,dotted] (-4,-2) ellipse (2 and 1.41);
\draw [blue,very thick,dotted] (-4,-2) ellipse (3.46 and 2.45);
\draw [blue,very thick,dotted] (-4,-2) ellipse (4.9 and 3.46);
\draw [blue,very thick,dotted] (-4,-2) ellipse (4.9 and 3.46) node [xshift=0.6cm, yshift=2.1cm] {level sets of $I(\xi)$};
\filldraw [fill=red!8,draw=none]
(4,-3) .. controls (2,0) and (1,-1)
.. (0,0).. controls (-1,1) and (-0.9,2) .. (-1,4) -- (4,4)  -- (4,-3);
\draw [blue,very thick,dotted] (-4,-2) ellipse (6.93 and 4.9);
\draw [blue,very thick,dotted] (-4,-2) ellipse (8.49 and 6);  	
\draw [orange, very thick] (3,-3) -- (-4,4) %
node[fill=white,near end, sloped, above, xshift=-5pt]{$F_1(u,\xi;\xistar)=z$}
;
\draw[red!60!black,very thick] (-1,4)  .. controls (-0.9,2)and (-1,1)
.. (0,0).. controls (1,-1) and (2,0) .. (4,-3) node[near end, sloped, below]
{$F(u,\xi)=z$};

\draw [red, line width=0.7mm] 	(4,-1) .. controls (3.95,-0.95) and (1,-1) .. (0,0).. controls (-1,1) and (-0.95,3.95) .. (-1,4) %
	node[fill=red!8,yshift=-1cm,xshift=1.7cm]{$F_2(u,\xi;\xistar)=z$}
	;  	
\draw[black, yshift=3.5cm,xshift=6cm] node[fill=red!8] {\textcolor{black}{$\{\xi: F(u,\xi)\geq z\}$}}; 		
\filldraw [black] (0,0) circle (2pt) node[yshift=-0.5cm,xshift=0cm]{$\xistar$};
\draw [black, very thick, ->] (0,0)	-- (1,1) node[yshift=-0.3cm,xshift=0.3cm] {$\nstar$};
\end{axis}
\end{tikzpicture}
\caption{Illustration of first and second-order approximations of probability using Taylor expansions at the LDT minimizer $\xistar$, which is the solution to \eqref{eq:LDT-min-opt} . The blue dotted lines are level sets of the rate function
	$I(\cdot)$, and the red shaded region is the extreme event set $\{\xi: F(u,\xi)\geq z\}$.  The normal $\nstar$ aligns with the gradients $\nabla_\xi F(u,\xistar)$
	and $\nabla I(\xistar)$.
	The figure is adapted from \cite{tong2020extreme}.
}\label{fig:LDT}
\end{figure}
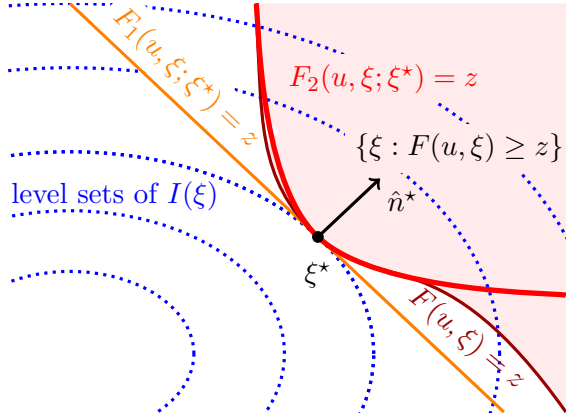

%% file: gaussian.tex
\section{Explicit formulations for Gaussian random parameter}\label{sec:G}

Throughout this section, we shall assume that $\xi\sim\Gauss{\mu}{\Sigma}$ follows a Gaussian distribution with mean $\mu$ and covariance $\Sigma \succ 0$.
Our goal will be to provide explicit formulations of the first- and second-order LDT chance constrained optimization problems \eqref{eq:op-LDT}.

\subsection{Rate function}
The moment generating function of $\xi$ is $\exp\left(\eta^\top \mu+\frac{1}{2}\eta^\top\Sigma\eta\right)$, and therefore, its cumulant generating function, $S(\eta) = \eta^\top \mu+\frac{1}{2}\eta^\top\Sigma\eta$, is quadratic.
Its rate function can be computed explicitly through the Legendre transform of $S$, and is also quadratic:
\begin{equation}
\label{eq:I-G}
I(\xi) = \frac{1}{2}\|\xi-\mu\|_{\SigInv{}}^2.
\end{equation}
Observe that $I(\xi)$ is, up to a normalization constant, the
negative log-probability density of $\Gauss{\mu}{\Sigma}$.
Hence, for a Gaussian distribution, %
minimizing the rate function over some domain, such as in \eqref{eq:LDT-min}, 
is equivalent to maximizing the probability density function over that domain.
We note, however, this is a special property of Gaussians, and it does not hold in general;
indeed, the next section shows that this property no longer holds for non-Gaussian parameters.

\subsection{First- and second-order probability estimates} \label{sec:prob-G}
For $\xi \sim\Gauss{\mu}{\Sigma}$, the probability estimates in \eqref{eq:PFOPSO} are Gaussian integrals over a half-space ($k = 1$), or a set bounded by a quadric ($k = 2$).
The explicit formulas for $P_k$ can be computed by first splitting the respective domains into a normal and orthogonal component (with respect to $\nabla_{\xi} F(u, \xistar)$), and then by computing the integral along each direction directly \cite{tong2020extreme}.
In the remainder of this discussion, we fix $u$ and $\xistar$ and assume that they satisfy $F(u, \xistar) = z$.

The first-order approximation is a Gaussian integral over a half-space bounded by a plane perpendicular to $\nstar \coloneqq \nabla_{\xi} F(u, \xistar)/\|\nabla_{\xi} F(u, \xistar)\|$, and passing through $\xistar$.
It is straightforward to verify that this integral degenerates to a one-dimensional Gaussian integral and can be computed using the dominating point $\xistar$ as follows:
\begin{equation}\label{eq:PFO-G}
P_1(u,z,\xistar)=\Phi(-\sqrt{2I(\xistar)}) = \Phi(-\|\xistar-\mu\|_{\SigInv{}}).
\end{equation}

The second-order probability estimate $P_2$ is slightly more complicated, but can still be approximated using an analytical expression. %
Following \cite{tong2020extreme}, we first approximate the surface $F_2(u,\xi; \xistar) = z$ by a paraboloid with axis of symmetry $\nstar$ at $\xistar$, and curvatures from $\nabla_\xi^2 F(u, \xistar)$ in the directions on the orthogonal plane.
$P_2(u,z,\xistar)$ can then be approximated as the measure of the set bounded by that paraboloid, which turns out to be equal to $P_1(u,z,\xistar)$ multiplied by a correction term. This correction term can be computed either through principle curvatures of $F$ at $\xistar$ \cite{du2001most}, or by using the eigenvalues of the rotated covariance-preconditioned Hessian of $F$ at $\xistar$; see \cite{tong2020extreme} for details.
In this paper, we introduce another way to compute this correction term, which turns out to be more numerically tractable for optimization purposes.
 \begin{theorem}[Second-order probability estimate for Gaussians]\label{thm:FO-G}
     Suppose that $\xi\sim\Gauss{\mu}{\Sigma}$, $u \in \mathcal{U}$ and $\xistar \in \Xi$ satisfy $F(u, \xistar) = z$, and that the following matrix is invertible:
     \[
     H \coloneqq \mathrm{I}_n-\frac{\|\xistar-\mu\|_{\SigInv{}} }{\left\|\SigHalf{}\nabla_\xi F(u,\xistar)\right\|}\SigHalf{} \nabla^2_\xi F(u,\xistar)\SigHalf{}
     \]
     Then the second-order approximation of the probability $\PP(F(u,\xi)\geq z)$ is given by:
 	 \begin{equation}\label{eq:PSO-G}
 	 P_2(u,z,\xistar) \approx
 	 \Phi(-\left\|\xistar-\mu\right\|_{\SigInv{}})
      \det\nolimits_{\perp \nhat} \left(H\right)^{-\frac{1}{2}},
 	 \end{equation}
      where the covariance-preconditioned normal $\nhat = \dfrac{\SigHalf{}\nabla_\xi F(u,\xistar)}{\|\SigHalf{}\nabla_\xi F(u,\xistar)\|} $ and orthogonal determinant %
 	 \begin{equation}\label{eq:detp}
 	 \det\nolimits_{\perp\nhat} H %
      = \left(\hat{n}^\top H^{-1} \hat{n}\right)\det H.
 	 \end{equation}
 \end{theorem}

\begin{proof}    
	Following the proof of Theorem~4.4 in~\cite{tong2020extreme}, the second-order probability approximation can be shown to be equal to:
	\begin{equation*}
	P_2(u,z,\xistar) = P_1(u,z,\xistar) C(u, \xistar),
	\end{equation*}
	where $C(u, \xistar)$ is a correction term computed through the curvatures of $F$ at $\xistar$.
    To define it, we first denote $S(\vartheta)$ to be the Lebesgue measure on the orthogonal component $\vartheta\in\left\lbrace \vartheta\in\RR^n: \<\nhat,\vartheta\> =0\right\rbrace$ of %
    $\nhat$ at $\xistar$.
    It can then be shown that
	\begin{equation*}
    \begin{aligned}
    C(u, \xistar) \approx& (2\pi)^{-\frac{n-1}{2}}\int_{\left\lbrace \vartheta\in\RR^n: \<\nhat,\vartheta\> =0\right\rbrace}
    \exp\left(
            -\frac{1}{2}\<\vartheta, \mathrm{I}_n - \frac{\|\SigHalf{}\nabla I(\xistar)\|}{\|\SigHalf{}\nabla_\xi F(u,\xistar)\|} \SigHalf{}\nabla_\xi^2 F(u,\xistar) \SigHalf{}   \vartheta \>
        \right) \mathrm{d}S(\vartheta)\\
    =& (2\pi)^{-\frac{n-1}{2}}\int_{\left\lbrace \vartheta\in\RR^n: \<\nhat,\vartheta\> =0\right\rbrace}
    \exp\left(
    -\frac{1}{2}\<\vartheta, H   \vartheta \>
    \right) \mathrm{d}S(\vartheta)
    \end{aligned}
	\end{equation*}
    where we have used $\nabla I(\xistar)=\SigInv{}(\xistar-\mu)$ along with the definition of $H$.
	The above expression can be further simplified by
	defining the unit vector:
	\begin{equation}\label{eq:hatt}
	\hat{t}:=\frac{H^{-1} \nhat}{\|H^{-1} \nhat\|}.
	\end{equation}
	Then, for all $\xi \in \RR^n$, there exist $s\in\RR$ and $\vartheta\in\left\lbrace \vartheta\in\RR^n: \<\nhat,\vartheta\> =0\right\rbrace$ such that $\xi = \vartheta+ s\hat{t}$. From Fubini's theorem and the orthogonality between $\nhat$ and $\xi$, we have
	\begin{equation}\label{eq:Gauss-int}
	\begin{aligned}
	\int_{\RR^n} e^{-\frac{1}{2}\<\xi,H \xi\>}\mathrm{d} \xi =& \int_{\left\lbrace \vartheta\in\RR^n: \<\nhat,\vartheta\> =0\right\rbrace}\int_{\RR} e^{-\frac{1}{2}\<\vartheta+s\hat{t}, H (\vartheta+s\hat{t})\>} \<\nhat,\hat{t}\> \mathrm{d}s\, \mathrm{d}S(\vartheta)\\
	=&  \<\nhat,\hat{t}\> \int_{\RR}e^{-\frac{1}{2}\<\hat{t}, H \hat{t}\>s^2} \mathrm{d}s \int_{\left\lbrace \vartheta\in\RR^n: \<\nhat,\vartheta\> =0\right\rbrace} e^{-\frac{1}{2}\<\vartheta,H\vartheta\>} \mathrm{d}S(\vartheta).
	\end{aligned}
	\end{equation}
	Applying the Gauss integral formula %
    to both sides of \eqref{eq:Gauss-int}, we obtain
	\begin{equation*}
	(\det H)^{-\frac{1}{2}} = (2\pi)^{-\frac{n-1}{2}} \<\nhat,\hat{t}\>\<{\hat{t}}, H \hat{t}\>^{-\frac{1}{2}}
    \int_{\left\lbrace \vartheta\in\RR^n: \<\nhat,\vartheta\> =0\right\rbrace} e^{-\frac{1}{2}\<\vartheta,H\vartheta\>} \mathrm{d}S(\vartheta).
	\end{equation*}
	Plugging in the definition of $\hat{t}$ in \eqref{eq:hatt}, we obtain
	\begin{equation*}
	\begin{aligned}
    (2\pi)^{-\frac{n-1}{2}}\int_{\left\lbrace \vartheta\in\RR^n: \<\nhat,\vartheta\> =0\right\rbrace} e^{-\frac{1}{2}\<\vartheta,H\vartheta\>} \mathrm{d}S(\vartheta)
    =&\<\nhat,\hat{t}\>^{-1}  \<{\hat{t}}, H \hat{t}\>^{\frac{1}{2}} (\det H)^{-\frac{1}{2}}\\
	=& \left(\frac{\<\nhat,H^{-1}\nhat\>^2}{\|H^{-1} \nhat\|^2} \frac{\|H^{-1} \nhat\|^2}{\<H^{-1}\nhat,H H^{-1} \nhat\>} \det H \right)^{-\frac{1}{2}}\\
	=& \left(\left(\hat{n}^\top H^{-1} \hat{n}\right)\det H \right)^{-\frac{1}{2}} \\
    =& \det\nolimits_{\perp \nhat}\left(H\right)^{-\frac{1}{2}}.
	\end{aligned}
	\end{equation*}
    This establishes the statement of the theorem.
\end{proof}
Note that the expression \eqref{eq:PSO-G} in \Cref{thm:FO-G} is not the same as $P_2$ in \eqref{eq:PFOPSO}, and only provides an approximation of the latter. Nevertheless, under the assumption that $\xistar$ is the LDT minimizer in the lower-level problem of \eqref{eq:op-bi}, it is asymptotically equal to \eqref{eq:PFOPSO} as $z\to \infty$.
Crucially, it has the added advantage of being computable in closed-form.

\subsection{LDT chance-constrained problems}
An application of the approximations in \eqref{eq:PFO-G} and \eqref{eq:PSO-G} to \eqref{eq:op-LDT} yield the following first- ($k = 1$) and second-order ($k=2$) LDT approximations of the chance-constrained problem \eqref{eq:ccp}, respectively:
\begin{equation}\label{eq:op-G}
\begin{aligned}
& \underset{u, \xistar, \lambda}{\text{minimize}}
& & J(u) \\
& \text{subject to}
&& u \in \mathcal{U}, \; \xistar \in \Xi, \; \lambda \in \RR_{+} \\
&&& P_k(u,z,\xistar) \leq \alpha, \\ 
&&& F(u,\xistar)=z,\\
&&& \SigInv{}(\xistar-\mu)  = \lambda \nabla_\xi F(u,\xistar).
\end{aligned}
\end{equation}
Specifically, we replace the first-order probability constraint $P_1(u,z,\xistar) \leq \alpha$ as follows:
\begin{equation}\label{eq:op-G-FO}
-\sqrt{2I(\xistar)}\leq \Phi^{-1}(\alpha).
\end{equation}
Similarly, we replace the second-order probability constraint $P_2(u,z,\xistar) \leq \alpha$ as follows:
\begin{equation}\label{eq:op-G-SO}
\Phi(-\sqrt{2I(\xistar)})\det\nolimits_{\perp \nhat}\left(\mathrm{I}_n-\lambda\SigHalf{} \nabla^2_\xi F(u,\xistar)\SigHalf{}\right)^{-\frac{1}{2}} \leq \alpha,
\end{equation}
where we have used the fact that 
\begin{equation}\label{eq:lambda}
\frac{\left\|\xistar-\mu\right\|_{\SigInv{}} }{\left\|\SigHalf{}\nabla_\xi F(u,\xistar)\right\|}
=\frac{\left\|\SigHalf{} \lambda \nabla_\xi F(u,\xistar)\right\|}{\left\|\SigHalf{}\nabla_\xi F(u,\xistar)\right\|}
=\lambda.
\end{equation}

We close this section with a few important remarks.
First, note that since all of the functions involved in the closed-form expressions \eqref{eq:op-G-FO} and \eqref{eq:op-G-SO} are fairly standard, we can use standard automatic differentiation routines to compute their function values and gradients.
Notably, this allows us to use off-the-shelf solvers for the solution of \eqref{eq:op-G}.

Second, the matrix $H = \mathrm{I}_n-\lambda\SigHalf{} \nabla^2_\xi F(u,\xistar)\SigHalf{}$ appearing in \eqref{eq:op-G-SO} is guaranteed to be positive definite (and hence, invertible) %
whenever the local solutions of the lower-level problem in \eqref{eq:op-bi} satisfy the sufficient second-order optimality conditions.
Although the solutions of \eqref{eq:op-G} cannot guarantee this in general, we have nonetheless empirically found in our numerical experiments that the linear system $Hx = \hat{n}$ was almost always solvable; %
in such cases, we can equivalently reformulate \eqref{eq:detp} as $\hat{n}^\top x \det(H)$ and evaluate its function values and gradients.

Third, and finally, 
although the above derivations are for Gaussian random parameters, the probability estimates \eqref{eq:PFO-G} and \eqref{eq:PSO-G} as well as the resulting LDT chance-constrained problem \eqref{eq:op-G} can be also used for random parameters that can be transformed to a Gaussian distribution. Specifically, if for a random parameter $\zeta$, there exists an invertible function $f$, such that $f(\zeta)$ is Gaussian, then we can use $\xi=f(\zeta)$ instead as the random parameter and apply the above results to $\xi$ while involving $f^{-1}$ in $F$. For example, when $\zeta$ follows a log-normal distribution, we can use $\xi=\log(\zeta)$ as the random parameter, while replacing $\zeta$ with $\exp(\xi)$ in $F$. We discuss these techniques further in \cref{sec:shortcolumn,sec:stock}.

%% file: gaussian_mixture.tex
\section{Explicit formulations for Gaussian mixture random parameter}\label{sec:GM}
We now generalize the approach from the last section for non-Gaussian random parameters to provide explicit formulations of the first- and second-order LDT chance constrained optimization problems \eqref{eq:op-LDT}.
Specifically, we assume that $\xi\sim\sum_{i=1}^M w_i\Gauss{\mu_i}{\Sigma_i}$ follows a Gaussian mixture distribution with $M$ components, where $w_i > 0$ and $\sum_{i=1}^M w_i=1$. Here, we use $\Pm_i = \Gauss{\mu_i}{\Sigma_i}$ to denote the distribution of the $i^\text{th}$ component with mean $\mu_i$ and covariance $\Sigma_i \succ 0$. Note that $\Pm$ is a convex combination of $\Pm_i$, $i \in \{1, \ldots, M\}$:
\begin{equation}
\label{eq:P-GM}
\Pm = \sum\limits_{i=1}^M w_i \Pm_i.
\end{equation}

\subsection{Rate function}
Unlike in the case of Gaussians, the rate function is no longer equal to the negative log-density of $\xi$, and is available only implicitly.
We first compute the cumulant generating function using its definition as follows:
\begin{equation}
\label{eq:S-GM}
\begin{aligned}
S(\eta)=& \log\left\lbrace\int_{\RR^n} e^{\eta^\top \xi} \mathrm{d}\Pm(\xi) \right\rbrace= \log\left\lbrace\int_{\RR^n} e^{\eta^\top \xi} \sum_{i=1}^M w_i \mathrm{d}\Pm_i(\xi) \right\rbrace\\
=& \log\left\lbrace \int_{\RR^n} e^{\eta^\top \xi} \sum_{i=1}^M \left[w_i(2\pi)^{-\frac{n}{2}} (\det \Sigma_i)^{-\frac{1}{2}} e^{-\frac{1}{2}(\xi-\mu_i)^\top \Sigma_i^{-1} (\xi-\mu_i)}\right] \mathrm{d}\xi\right\rbrace \\
=&\log \left\lbrace \sum_{i=1}^M w_i e^{\eta^\top \mu_i+\frac{1}{2}\eta^\top\Sigma_i\eta  } \right\rbrace.
\end{aligned}
\end{equation}
Since the rate function $I$ is the Legendre transform of $S$, it is non-trivial to compute it explicitly in this case.
Nevertheless, as noted in \Cref{sec:LDT}, we can still compute the LDT minimizer $\xistar$ by solving \eqref{eq:xistar-eta} instead:
\begin{equation}
\label{eq:xistar-GM}
\begin{aligned}
(\xistar, \etastar) \in &\argmin_{\xi, \eta \in \RR^n} && \eta^\top\xi -\log \left\lbrace \sum_{i=1}^M w_i e^{\eta^\top \mu_i+\frac{1}{2}\eta^\top\Sigma_i\eta  } \right\rbrace\\
& \text{subject to} && F(u,\xi) = z, \quad \xi = \frac{\sum\limits_{i=1}^M w_i e^{{\eta}^\top\mu_i+\frac{1}{2}{\eta}^\top\Sigma_i\eta}\left(\mu_i+\Sigma_i \eta\right)  }{\sum\limits_{i=1}^M w_i e^{{\eta}^\top\mu_i+\frac{1}{2}{\eta}^\top\Sigma_i\eta}  }.
\end{aligned}
\end{equation}

Moreover, even though the rate function is not available in closed form, the following theorem establishes that Gaussian mixture distributions (and hence, Gaussian distributions) always satisfy condition \ref{assume:finte_volume} in \Cref{thm:ldt_verifiable}.
Therefore, the conditions of \Cref{thm:ldt_verifiable} depend only on the behavior of the constraint function $F$ with respect to the uncertain parameters $\xi$, and not on the form of the mixture distribution.

\begin{theorem}[Sufficient conditions for Gaussian mixture distributions]\label{thm:gaussian_mixture_finite_volume}
    Suppose that
    $\xi \sim \sum_{j = 1}^M w_i \Gauss{\mu_j}{\Sigma_j}$. %
    Then condition \ref{assume:finte_volume} in \Cref{thm:ldt_verifiable} is satisfied.
\end{theorem}
\begin{proof}
    \begingroup
    \allowdisplaybreaks
    We shall show that
    \begin{enumerate}[leftmargin=*,label=(\alph*)]
        \item\label{toshow:strongconvex} the cumulant generating function $S$ is strongly convex, and
        \item\label{toshow:ratelowbound} the rate function $I$ satisfies $I(\xi) \geq \min_{j \in \{1, \ldots, M\}} I_j(\xi)$ for all $\xi \in \RR^n$, where $I_j(\xi) = \frac{1}{2}\|\xi-\mu_j\|_{\SigInv{j}}^2$ denotes the rate function of the $j^\text{th}$ Gaussian $\Gauss{\mu_j}{\Sigma_j}$, see \cref{eq:I-G}.
    \end{enumerate}
    Assume for the moment that the above statements are true, and let $\sigma > 0$ denote the strong convexity parameter of $S$ in \ref{toshow:strongconvex}.
    Since the rate function~\cref{eq:I} is the convex conjugate of $S$,
    strong convexity of $S$ implies that the gradient of the rate function $\nabla I(\xi)$ is Lipschitz continuous with constant $\sigma$; that is,
    \begin{equation}\label{eq:lipschitz}%
    I(\hat{\xi}) + 
    \<\nabla I(\hat{\xi}), \xi-\hat{\xi}\>
    \geq I(\xi) - \frac{\sigma}{2}\|\xi - \hat{\xi}\|^2, \quad \forall \xi, \hat{\xi} \in \RR^n.
    \end{equation}
    Now, fix $r > 0$ and $\hat{\xi} \in \RR^n$, and let $\txi$ be the (unique) minimizer of the integrand in \ref{assume:finte_volume} over $B(r, \hat{\xi})$. Since for the given Gaussian mixture distribution $\Pm$, (see \Cref{sec:GM})
    \[
    \mathrm{d}\Pm(\xi) = \sum_{j=1}^M \left[w_j(2\pi)^{-\frac{n}{2}} (\det \Sigma_j)^{-\frac{1}{2}} \exp\left(-\frac{1}{2}\|\xi-\mu_j\|_{\SigInv{j}}^2 \right)\right] \mathrm{d}\xi,
    \]
    we have that
    \[
    \txi = \argmin_{\xi \in B(r, \hat{\xi})} \left\{
    \exp\left(
    I(\hat{\xi}) + 
    \<\nabla I(\hat{\xi}), \xi-\hat{\xi}\> \right)
    \sum_{j=1}^M \left[w_j(2\pi)^{-\frac{n}{2}} (\det \Sigma_j)^{-\frac{1}{2}} \exp\left(-\frac{1}{2}\|\xi-\mu_j\|_{\SigInv{j}}^2 \right)\right] %
    \right\}.
    \]
    Combining the above with \cref{eq:lipschitz}, we obtain
    \begin{alignat*}{2}
    \int_{\mathcal{B}(r, \hat{\xi})} \exp\left(
    I(\hat{\xi}) + 
    \<\nabla I(\hat{\xi}), \xi-\hat{\xi}\> \right)
    \mathrm{d}\Pm(\xi)
    &\geq    
    \sum_{j = 1}^M c_j \exp\left(
    I(\txi) - I_j(\txi) - \frac{\sigma}{2}\|\txi - \hat{\xi}\|^2
    \right) \\
    &\geq
    \sum_{j = 1}^M c_j e^{-\frac{\sigma}{2}r^2} \exp\left(
    I(\txi) - I_j(\txi)
    \right), \\
    \Big(\text{with } c_j &\coloneqq
    w_j(2\pi)^{-\frac{n}{2}} (\det \Sigma_j)^{-\frac{1}{2}} 
    \int_{\mathcal{B}(r, \hat{\xi})} \mathrm{d}\xi > 0\Big),
    \end{alignat*}
    where we have used the definition of the $j^\text{th}$ component rate function  $I_j(\xi) = \frac{1}{2}\|\xi-\mu_j\|_{\SigInv{j}}^2$, along with the fact that $\txi, \hat{\xi} \in B(r, \hat{\xi}) \implies \|\txi - \hat{\xi}\| \leq r$.
    Claim \ref{toshow:ratelowbound} implies that there exists some $k \in \{1, \ldots, M\}$ such that $I(\txi) \geq I_k(\txi)$.
    Therefore, we have shown that
    \begin{alignat*}{2}
    \int_{\mathcal{B}(r, \hat{\xi})} \exp\left(
    I(\hat{\xi}) + 
    \<\nabla_\xi I(\hat{\xi}), \xi-\hat{\xi}\> \right)
    \mathrm{d}\Pm(\xi)
    &\geq 
    \min_{j \in \{1, \ldots, M\}} c_j e^{-\frac{\sigma}{2}r^2} > 0, \quad \forall \hat{\xi} \in \RR^n.
    \end{alignat*}
    Also, we have
    \begin{equation*}
    \int_{\mathcal{B}(r, \hat{\xi})} \exp\left(
    I(\hat{\xi}) + 
    \<\nabla I(\hat{\xi}), \xi-\hat{\xi}\> \right)
    \mathrm{d}\Pm(\xi)
    \leq 
    \int_{\RR^n} \exp\left(
    I(\hat{\xi}) + 
    \<\nabla I(\hat{\xi}), \xi-\hat{\xi}\> \right)
    \mathrm{d}\Pm(\xi)
    = 1.
    \end{equation*}
    where the equality can be verified easily (\textit{e.g.}, see the proof of Theorem~2.1~in~\cite{dematteis2019extreme}).
    Therefore, the left-hand side expression in condition \ref{assume:finte_volume} is bounded away from $0$ by a constant that is independent of $\hat{\xi}$, whereas its right-hand side $\to \infty$ as $\|\hat{\xi}\| \to \infty$ (see proof of \Cref{thm:ldt_verifiable}).
    This proves the statement of this theorem.
    
    It remains to prove claims \ref{toshow:strongconvex} and \ref{toshow:ratelowbound}.
    To show \ref{toshow:strongconvex}, let $S_j(\eta) = \eta^\top \mu_j+\frac{1}{2}\eta^\top\Sigma_j\eta$ denote the cumulant generating function of the $j^\text{th}$ Gaussian $\Gauss{\mu_j}{\Sigma_j}$.
    Then, the cumulant generating function of the Gaussian mixture $S(\eta) = \log \left(\sum_{j = 1}^M w_j \exp(S_j(\eta)) \right)$, see \cref{eq:S-GM}.
    A straightforward calculation shows that its Hessian satisfies
    \begin{alignat*}{2}
    \nabla^2 S(\eta) &= \left({
        \begin{aligned}
        \left[\sum_{j=1}^M e^{S_j(\eta)}\right] \left[\sum_{j=1}^M e^{S_j(\eta)} \left(\nabla S_j(\eta) (\nabla S_j(\eta))^\top + \nabla^2 S_j(\eta) \right)\right] \\
        -
        \left[\sum_{j=1}^M e^{S_j(\eta)} \nabla S_j(\eta) \right] \left[\sum_{j=1}^M e^{S_j(\eta)} (\nabla S_j(\eta))^\top \right]
        \end{aligned}
    }\right)\cdot{
        \left[\sum_{j=1}^M e^{S_j(\eta)}\right]^{-2}
    } \\
    &=
    \left({
        \begin{aligned}
        \sum_{j=1}^M e^{2S_j(\eta)} \nabla^2 S_j(\eta) + \sum_{j=1}^M\sum_{i=j+1}^M e^{S_i(\eta) + S_j(\eta)} \Big(
        \nabla^2 S_i(\eta) + \nabla^2 S_j(\eta) \\+ \left[\nabla S_i(\eta) - \nabla S_j(\eta) \right] \left[\nabla S_i(\eta) - \nabla S_j(\eta) \right]^\top
        \Big)
        \end{aligned}
    }\right)\cdot{
        \left[\sum_{j=1}^M e^{S_j(\eta)}\right]^{-2}
    } \\
    &\succeq
    \varepsilon \mathrm{I}_n \cdot \left({
        \begin{aligned}
        \sum_{j=1}^M e^{2S_j(\eta)} + \sum_{j=1}^M\sum_{i=j+1}^M 2 e^{S_i(\eta) + S_j(\eta)} 
        \end{aligned}
    }\right)\cdot{
        \left[\sum_{j=1}^M e^{S_j(\eta)}\right]^{-2}
    } \\
    & = \varepsilon \mathrm{I}_n,
    \end{alignat*}
    where $\varepsilon = \min_{j \in \{1, \ldots, M\}} \lambda_{\text{min}} (\nabla^2 S_j) =  \min_{j \in \{1, \ldots, M\}} \lambda_{\text{min}} (\Sigma_j) > 0$, and $\lambda_{\text{min}}(A)$ denotes the minimum eigenvalue of a matrix $A$.
    
    To show \ref{toshow:ratelowbound}, let $f^*$ denote the conjugate of a function $f$.
    Convexity and continuity of the cumulant generating function implies $I^{*} = (S^*)^{*} = I$,
    and hence, for any $\eta \in \RR^n$, we have
    \begin{alignat*}{2}
    I^*(\eta) =
    S(\eta) = \log\left( \sum_{j = 1}^M w_j \exp(S_j(\eta)) \right)
    &\leq \log\left( \max_{j \in \{1, \ldots, M\}} \exp(S_j(\eta)) \right) \\
    &= \max_{j \in \{1, \ldots, M\}} S_j(\eta) \\
    &= \max_{j \in \{1, \ldots, M\}} \max_{\xi \in \RR^n} \left\{ \eta^\top \xi - I_j(\xi) \right\} \\
    &= \max_{\xi \in \RR^n} \left\{ \eta^\top \xi - \min_{j \in \{1, \ldots, M\}} I_j(\xi) \right\} \\
    &= \left(\min_{j \in \{1, \ldots, M\}} I_j \right)^*(\eta)
    \end{alignat*}
    where the first inequality follows from the fact that $w_j \geq 0$, $\sum_{j=1}^M w_j = 1$, the first equality follows from the monotonicity of $\exp$, the second and final equalities follow from the definition of the convex conjugate, 
    and the third equality follows because we can interchange the order of the $\max$ over $j$ with the $\max$ over $\xi$.
    Finally, the order reversing property of the conjugate $f^* \leq g^* \implies f \geq g$, implies that $I(\xi) \geq \min_{j \in \{1, \ldots, M\}} I_j (\xi)$ which is claim \ref{toshow:ratelowbound}.
    \endgroup
\end{proof}

\subsection{First- and second-order probability estimates}\label{sec:prob-GM}
Recall from \cref{sec:prob} that the first-order probability estimate $P_1$ in \eqref{eq:PFOPSO} is the measure of the half-space bounded by $F_1(u,\xi; \xistar) = z$, where $\xistar$ solves \eqref{eq:xistar-GM}.
Since a Gaussian mixture is just a convex combination of Gaussians, we can exploit the analogous result for the latter~\eqref{eq:PFO-G} to compute the first-order estimate analytically.

\input{GMM_Illustration.tex}

\begin{theorem}[First-order probability estimate for Gaussian mixtures]\label{thm:FO-GM}
	Suppose that $\xi\sim\sum_{i=1}^M w_i\Gauss{\mu_i}{\Sigma_i}$ and for some fixed $u \in \mathcal{U}$, $\xistar$ is the solution to \eqref{eq:xistar-GM}. Then the first-order approximation of the probability $\PP(F(u,\xi)\geq z)$ is given by
	\begin{equation}\label{eq:PFO-GM}
	P_1(u,z,\xistar)= \sum\limits_{i=1}^M w_i \Phi\left(-\frac{\<\nabla_\xi F(u, \xistar),\xistar-\mu_i\>}{\|\SigHalf{i}\nabla_\xi F(u,\xistar )\|}\right).
	\end{equation}
\end{theorem}

\begin{proof}
    \begingroup
    \allowdisplaybreaks
Substituting \eqref{eq:P-GM} in \eqref{eq:PFOPSO} for $k=1$, we obtain:
\begin{equation}\label{eq:PFO-linear}
P_1(u,z,\xistar)
= \sum_{i=1}^{M} w_i \Pm_i\left(\left\{\xi: F_1(u,\xi; \xistar) \geq z \right\}\right).
\end{equation}
Therefore, we can consider each $\Pm_i\left(\left\{\xi: F_1(u,\xi; \xistar) \geq z \right\}\right)$, $i \in \{1, \ldots, M\}$, separately.
Similar to \eqref{eq:PFO-G} for Gaussians, this quantity is a Gaussian integral over a half-space.
However, unlike the latter,
the contours of negative log-density of $\Gauss{\mu_i}{\Sigma_i}$, \textit{i.e.}, of $\frac{1}{2}\|\xi-\mu_i\|_{\SigInv{i}}^2$, need not necessarily be tangential to the half-space $F_1(u,\xi; \xistar) = z$ at $\xistar$.
Therefore, to compute $\Pm_i\left(\left\{\xi: F_1(u,\xi; \xistar) \geq z \right\}\right)$, we need to find the point $\txi_i$ on the hyperplane $F_1(u,\xi; \xistar)=z$ where the contour of $\frac{1}{2}\|\xi-\mu_i\|_{\SigInv{i}}^2$ is tangential to the plane (see~\Cref{fig:GM}):
\begin{equation}
\label{eq:txi1}
\begin{aligned}
\txi_i
&= \argmin_{\xi \in \RR^n} \left\{\frac{1}{2}\|\xi-\mu_i\|_{\SigInv{i}}^2 : F_1(u,\xi; \xistar)=z\right\} \\
&= \argmin_{\xi \in \RR^n} \left\{\frac{1}{2}\|\xi-\mu_i\|_{\SigInv{i}}^2 : \<\nabla_\xi F(u,\xistar), \xi-\xistar \>=0\right\},
\end{aligned}
\end{equation}
where we used the definition of the first-order Taylor expansion \eqref{eq:FFO} and the fact that $\xistar$ satisfies $F(u,\xistar) = z$ since it solves \eqref{eq:xistar-GM}.

The optimality conditions of \eqref{eq:txi1} imply that there exists $\lambda_i \geq 0$ such that: 
\begin{equation*}
\SigInv{i}(\txi_i-\mu_i) = \lambda_i \nabla_\xi F(u,\xistar).
\end{equation*}
Therefore, $\txi_i=\mu_i+\lambda_i \Sigma_i \nabla_\xi F(u,\xistar)$.
Substituting this in $\<\nabla_\xi F(u,\xistar), \txi_i-\xistar \>=0$, we obtain
\[
\lambda_i=\frac{\< \nabla_\xi F(u,\xistar), \xistar - \mu_i\>}{\nabla_\xi F(u,\xistar), \Sigma_i \nabla_\xi F(u,\xistar)}.
\]
Therefore, we have that
\[
\txi_i = \mu_i + \frac{\< \nabla_\xi F(u,\xistar), \xistar - \mu_i\>}{\<\nabla_\xi F(u,\xistar), \Sigma_i \nabla_\xi F(u,\xistar)\>} \Sigma_i \nabla_\xi F(u,\xistar).
\]
Finally, replacing $\xistar$ in the formula of the first-order probability estimate for Gaussians \eqref{eq:PFO-G} with $\txi_i$ yields
\begin{equation*}
\begin{aligned}
\Pm_i\left(\left\{\xi: F_1(u,\xi; \xistar) \geq z \right\}\right) 
&= \Phi\left(-\left\|\txi_i-\mu_i\right\|_{\SigInv{i}}\right)\\
&= \Phi\left(-\left\| \mu_i + \frac{\< \nabla_\xi F(u,\xistar), \xistar - \mu_i\>}{\<\nabla_\xi F(u,\xistar), \Sigma_i \nabla_\xi F(u,\xistar)\>} \Sigma_i \nabla_\xi F(u,\xistar)-\mu_i\right\|_{\SigInv{i}}\right)\\
&= \Phi\left(- \frac{\< \nabla_\xi F(u,\xistar), \xistar - \mu_i\>}{\<\nabla_\xi F(u,\xistar), \Sigma_i \nabla_\xi F(u,\xistar)\>}\left\| \Sigma_i \nabla_\xi F(u,\xistar) \right\|_{\SigInv{i}}\right)\\
&= \Phi\left(-\frac{\<\nabla_\xi F(u, \xistar),\xistar-\mu_i\>}{\|\SigHalf{i}\nabla_\xi F(u,\xistar )\|}\right).
\end{aligned}
\end{equation*}
Combining this with \eqref{eq:PFO-linear} gives us the claimed formula \eqref{eq:PFO-GM}.
\endgroup
\end{proof}

The second-order probability estimate $P_2$ satisfies a linear relationship similar to the first-order estimate $P_1$; that is,
$
P_2(u,z,\xistar)
= \sum_{i=1}^{M} w_i \Pm_i\left(\left\{\xi: F_2(u,\xi; \xistar) \geq z \right\}\right).
$
For each $i \in \{1, \ldots, M\}$, we can then use formula \eqref{eq:PSO-G} for Gaussians to compute $\Pm_i\left(\left\{\xi: F_2(u,\xi; \xistar) \geq z \right\}\right)$.
As in the case of the first-order estimate, this requires us to first find the point on the surface $F_2(u,\xi; \xistar)=z$ where the contours of $\frac{1}{2}\|\xi-\mu_i\|_{\SigInv{i}}^2$ are tangential to the surface, see \Cref{fig:GM}.
The final closed-form expression is presented in the following theorem, whose proof is omitted for brevity.

\begin{theorem}[Second-order approximation of probability for Gaussian mixtures]\label{thm:SO-GM}
    Suppose that $\xi\sim\sum_{i=1}^M w_i\Gauss{\mu_i}{\Sigma_i}$ and for some fixed $u \in \mathcal{U}$, $\xistar$ is the solution to \eqref{eq:xistar-GM}. Then the second-order approximation of the probability $\PP(F(u,\xi)\geq z)$ is given by
	\begin{equation}\label{eq:PSO-GM}
	P_2(u,z,\xistar)=\sum\limits_{i=1}^M w_i \Phi\left(-\|\txi_i-\mu_i\|_{\SigInv{i}}\right)
	\det\nolimits_{\perp \tn_i}\left(\mathrm{I}_n-\frac{\|\txi_i-\mu_i\|_{\SigInv{i}} }{\|\SigHalf{i}\nabla_\xi F_2(u,\txi_i; \xistar)\|}\SigHalf{i} \nabla^2_\xi F(u,\xistar)\SigHalf{i}\right)^{-\frac{1}{2}},
	\end{equation}
	where %
	\begin{equation}\label{eq:txi2}
	\txi_i= \argmin_{\xi \in \RR^n} \left\{\frac{1}{2}\|\xi-\mu_i\|_{\SigInv{i}}^2 : F_2(u,\xi; \xistar)=z\right\}, \quad \tn_i = \frac{\SigHalf{i} \nabla_\xi F_2(u, \txi_i; \xistar)}{\|\SigHalf{i}\nabla_\xi F_2(u, \txi_i; \xistar)\|}.
	\end{equation}
\end{theorem}

\subsection{LDT chance-constrained problems}
Substituting \eqref{eq:xistar-GM} and \eqref{eq:PFO-GM} in \eqref{eq:op-LDT-etastar} ($k = 1$) yields the following first-order LDT approximation of the chance-constrained problem \eqref{eq:ccp}:
\begin{equation}\label{eq:op-GM1}
\begin{aligned}
& \underset{u, \xistar, \etastar, \lambda}{\text{minimize}}
& & J(u) \\
& \text{subject to}
&& u \in \mathcal{U}, \; \xistar \in \RR^n, \; \etastar \in \RR^n, \; \lambda \in \RR_{+} \\
&&& \sum\limits_{i=1}^M w_i \Phi\left(-\frac{\<\nabla_\xi F(u, \xistar),\xistar-\mu_i\>}{\|\SigHalf{i}\nabla_\xi F(u,\xistar )\|}\right) \leq \alpha, \\
&&& F(u,\xistar)=z,\\
&&& \etastar  = \lambda \nabla_\xi F(u,\xistar), \\
&&& \xistar = \frac{\sum\limits_{i=1}^M w_i \exp\left({{\etastar}^\top\mu_i+\frac{1}{2}{\etastar}^\top\Sigma_i\etastar}\right)\left(\mu_i+\Sigma_i \etastar\right)  }{\sum\limits_{i=1}^M w_i \exp\left({{\etastar}^\top\mu_i+\frac{1}{2}{\etastar}^\top\Sigma_i\etastar}\right)  }.
\end{aligned}
\end{equation}

To obtain the second-order LDT approximation, we note that we first have to reformulate the nonconvex quadratic optimization problem \eqref{eq:txi2} that appears in the closed-form expression of the second-order probability estimate \eqref{eq:PSO-GM}.
The following lemma provides a characterization of $\txi_i$ that we shall use subsequently.
\begin{lemma}\label{lemma:qcqp}
Fix $u \in \mathcal{U}_0$ and suppose that $\xistar$ is the solution to \eqref{eq:xistar-GM}.
Then, for each $i \in \{1, \ldots, M\}$, $\txi_i$ is a global minimizer of
$
\min_{\xi \in \RR^n} \left\{\frac{1}{2}\|\xi-\mu_i\|_{\SigInv{i}}^2 : F_2(u,\xi; \xistar)=z\right\}
$
if and only if there exists $\tilde{\lambda}_i \in \RR_{+}$ such that
\begin{subequations}
    \begin{alignat}{1}
    &F_2(u,\txi_i; \xistar)=z, \label{eq:qcqp_1} \\
    &\SigInv{i}(\txi_i-\mu_i)=\tilde{\lambda}_i\left(\nabla_\xi F(u,\xistar)+\nabla_\xi^2 F(u,\xistar)(\txi_i-\xistar) \right), \label{eq:qcqp_2}\\
    &\mathrm{I}_n-\tilde{\lambda}_i\SigHalf{i} \nabla^2_\xi F(u,\xistar)\SigHalf{i} \succeq 0. \label{eq:qcqp_3}
    \end{alignat}
\end{subequations}
\end{lemma}
\begin{proof}
    First, note that $\nabla_{\xi} F_2(u,\xistar; \xistar) = \nabla_{\xi} F(u,\xistar) \neq 0$ (see \cref{sec:notation}).
    Therefore, \cite[Theorem~3.2]{more1993generalizations} shows that $\txi_i$ is a global minimizer if and only if \eqref{eq:qcqp_1} and \eqref{eq:qcqp_2} are satisfied, along with
    $
    \SigInv{i} - \tilde{\lambda}_i \nabla^2_\xi F_2(u, \txi_i; \xistar) \succeq 0.
    $
    Since $\nabla^2_\xi F_2(u, \txi_i; \xistar) = \nabla^2_\xi F(u, \xistar)$ by construction, this condition is equivalent to
    $
    \SigInvHalf{i}\left(\mathrm{I}_n-\tilde{\lambda}_i\SigHalf{i} \nabla^2_\xi F(u,\xistar)\SigHalf{i} \right)\SigInvHalf{i} \succeq 0,
    $
    which is seen to be equivalent to \eqref{eq:qcqp_3}.
\end{proof}

To obtain the second-order LDT chance-constrained approximation, we relax the positive definiteness condition \eqref{eq:qcqp_3}, and combine \eqref{eq:xistar-GM}, \eqref{eq:PSO-GM}, \eqref{eq:qcqp_1} and \eqref{eq:qcqp_2} along with \eqref{eq:op-LDT-etastar} ($k = 2$) to obtain the following formulation:
\begin{equation}\label{eq:op-GM2}
\begin{aligned}
& \underset{u, \xistar, \etastar, \lambda, \txi, \tilde{\lambda}}{\text{minimize}}
& & J(u) \\
& \text{subject to}
&& u \in \mathcal{U}, \; \xistar \in \RR^n, \; \etastar \in \RR^n, \; \lambda \in \RR_{+}, \\
&&& \txi_i \in \RR^n, \tilde{\lambda}_i \in \RR_{+}, \quad i \in \{1, \ldots, M\}, \\
&&& \<\nabla_\xi F(u,\xistar), \txi_i-\xistar \> +\frac{1}{2}\<\txi_i-\xistar, \nabla^2_\xi F(u,\xistar)( \txi_i-\xistar) \>=0, \quad i \in \{1, \ldots, M\}, \\
&&& \SigInv{i}(\txi_i-\mu_i)=\tilde{\lambda}_i\left(\nabla_\xi F(u,\xistar)+\nabla_\xi^2 F(u,\xistar)(\txi_i-\xistar) \right), \quad i \in \{1, \ldots, M\}, \\
&&& \sum\limits_{i=1}^M w_i \Phi\left(-\|\txi_i-\mu_i\|_{\SigInv{i}}\right)
\det\nolimits_{\perp \tn_i}\left(\mathrm{I}_n-\tilde{\lambda}_i\SigHalf{i} \nabla^2_\xi F(u,\xistar)\SigHalf{i}\right)^{-\frac{1}{2}}\leq \alpha,\\
&&& F(u,\xistar)=z,\\
&&& \etastar  = \lambda \nabla_\xi F(u,\xistar), \\
&&& \xistar = \frac{\sum\limits_{i=1}^M w_i \exp\left({{\etastar}^\top\mu_i+\frac{1}{2}{\etastar}^\top\Sigma_i\etastar}\right)\left(\mu_i+\Sigma_i \etastar\right)  }{\sum\limits_{i=1}^M w_i \exp\left({{\etastar}^\top\mu_i+\frac{1}{2}{\etastar}^\top\Sigma_i\etastar}\right)  }.
\end{aligned}
\end{equation}

As in the case of the LDT chance-constrained formulations for Gaussians, the function values and gradients of all functions involved in \eqref{eq:op-GM1} and \eqref{eq:op-GM2} can be computed using standard automatic differentiation codes, enabling their solution using off-the-shelf solvers.

%% file: GMM_Illustration.tex
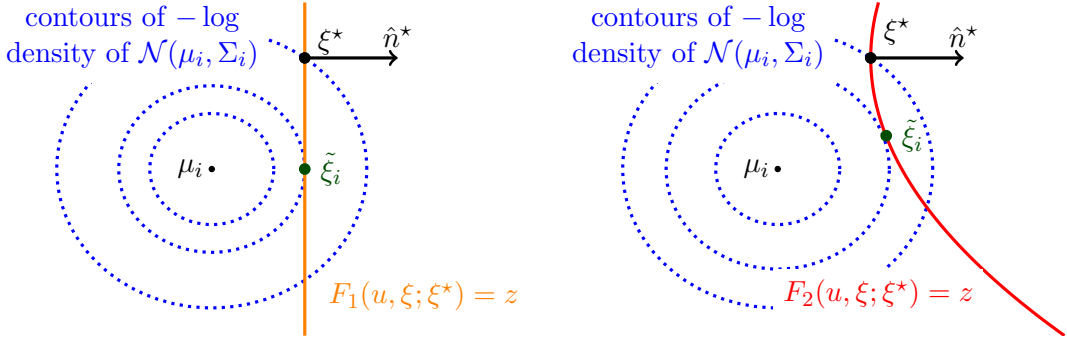
\begin{figure}[tbhp]\centering
	\begin{tikzpicture}[scale=1.]
	\begin{groupplot}[group style = {columns=2, rows=1, vertical sep =70pt, horizontal sep = 3pt}, width = 9cm, height = 6cm,compat=1.11,
	xmin=-7,
	xmax=11,
	ymin=-6,
	ymax=6,
	axis line style={draw=none},
	tick style={draw=none},
	yticklabels={,,},
	xticklabels={,,},
	]
	
	\nextgroupplot
	\draw  [blue,very thick,dotted]  (0,0) circle (2);
	\draw  [blue,very thick,dotted]  (0,0) circle (3);
	\draw  [blue,very thick,dotted]  (0,0) circle (5) node [fill=white,xshift=-30pt,yshift=50pt] 
	{\makecell[c]{contours of   $-\log$\\ density of $\mathcal{N}(\mu_i,\Sigma_i)$}};
	\draw [orange,very thick] (3,6) -- (3,-6) node[xshift=45pt,yshift=15pt] {$F_1(u,\xi; \xistar)=z$};
	\draw [->,very thick] (3,4) -- (6,4) node[above] {$\nstar$};
	\filldraw [black] (3,4) circle (2pt) node[xshift=10pt,yshift=6pt]{$\xistar$};
	\filldraw [green!30!black] (3,0) circle (2pt) 
	node[xshift=10pt]{$\txi_{i}$}
	;
	\filldraw [black] (0,0) circle (1pt) node[xshift=-8pt]{$\mu_i$};
	\nextgroupplot
	\draw  [blue,very thick,dotted]  (0,0) circle (2);
	\draw  [blue,very thick,dotted]  (0,0) circle (3.6);
	\draw  [blue,very thick,dotted]  (0,0) circle (5) node [fill=white,xshift=-30pt,yshift=50pt] 
	{\makecell[c]{contours of  $-\log$  \\ density of $\mathcal{N}(\mu_i,\Sigma_i)$}};
	\draw[red, very thick]   plot[smooth,domain=-6:6] ({(\x-4)^2/16+3}, \x) node[fill=white,yshift=-110pt,xshift=0pt] {$F_2(u,\xi; \xistar)=z$};
	\draw [->,very thick] (3,4) -- (6,4) node[above] {$\nstar$};
	\filldraw [black] (3,4) circle (2pt) node[xshift=10pt,yshift=10pt]{$\xistar$};
	\filldraw [green!30!black] (3.5,1.2) circle (2pt) 
	node[xshift=10pt]{$\txi_{i}$}
	;
	\filldraw [black] (0,0) circle (1pt) node[xshift=-8pt]{$\mu_i$};		
	\end{groupplot}
	\end{tikzpicture}
	\caption{\label{fig:GM}Illustration of finding $\txi_i$ for each Gaussian component $\Gauss{\mu_i}{\Sigma_i}$. The left figure shows the procedure for the first-order approximation, namely of finding a point on the hyperplane $F_1(u,\xi; \xistar)=z$ where the contour of negative log-density of $\Gauss{\mu_i}{\Sigma_i}$ is tangential to the plane, by solving \eqref{eq:txi1}. Similarly, the right figure shows the procedure for the second-order approximation, by finding $\txi_i$ on $F_2(u,\xi; \xistar)=z$ by solving \eqref{eq:txi2}.}
\end{figure}

%% file: applications.tex
\section{Applications} \label{sec:app}

We illustrate the computational performance of our first- and second-order LDT probability estimates and corresponding chance-constrained optimization problems using three examples from different applications.
In \cref{sec:stock}, we consider a portfolio selection problem with random stock prices;
in \cref{sec:shortcolumn}, we consider a structural short column design with uncertain material properties and loads;
and finally in \cref{sec:PDE}, we consider optimal control of a steady-state advection-diffusion PDE.
We implemented all of our formulations in Julia using JuMP~\cite{dunning2017jump}, and solved all optimization problems using Ipopt~\cite{wachter2006implementation} with default solver options.
We used the forward mode automatic differentiation provided by JuMP for gradient computations.
Since JuMP does not support Hessian computations of user-defined nonlinear constraint functions, which is needed for the second-order LDT chance-constrained formulations (for example, to encode the orthogonal determinant \eqref{eq:detp}),
we used Ipopt's in-built quasi-newton method for the solution of these problems.
It is important to note that our goal is to demonstrate the advantage and competitive performance of our method using a purely off-the-shelf implementation.
Reported solution times correspond to wall clock times on a single thread of an Intel~Xeon~Gold CPU at 2.30GHz with 512 GB RAM.

\subsection{Portfolio management}\label{sec:stock}
Consider an optimal asset allocation problem, where we wish to allocate our wealth among $n$ stocks which will be held over $T$ time periods.
Let $u_i \in \RR^{+}$ denote the fraction of wealth invested in the $i^\text{th}$ stock.
Our goal is to determine the optimal allocation $u=[u_1,\ldots,u_n]\in\RR^n_{+}$ of our wealth for each stock.
Let $v_i$ denote the value of the $i^\text{th}$ stock at the end of $T$ time periods.
We model the uncertainty in stock price $v_i$ using a traditional log-normal model \cite{hull2009options,kawas2011log},
\begin{equation}
\label{eq:s-Si}
v_i(\xi_i) = \exp\left(\left(\theta_i - \frac{\sigma_i^2}{2}\right)T + \sqrt{T} \xi_i\right), \quad i \in \{1, \ldots, n\},
\end{equation}
where $\theta_i\in\RR$ and $\sigma_i\in\RR$ represent the drift and infinitesimal standard deviation of the L{\'e}vy process for the $i^\text{th}$ stock, respectively, and $\xi_i \in \RR$ is a zero-mean continuously compounded rate of return which represents the true uncertainty driver.
Since choosing the right distribution for $\xi$ is an unresolved question (e.g., see~\cite{kawas2011log}), we model $\xi$ using Gaussian (as in traditional models) as well as Gaussian mixture distributions.
Our goal is to ensure, with very high probability, that the total return after $T$ time periods, $\sum_{i = 1}^n v_i(\xi_i) u_i$, is higher than some value $z$;
note that this is equivalent to constraining the probability $\PP(\sum_{i = 1}^n v_i (\xi_i) u_i \leq z)$.

\subsubsection{Probability estimation}%

Before we consider the chance-contrained formulation,
we empirically assess the accuracy of the first- and second-order LDT estimates of the rare event probability $\PP(\sum_{i = 1}^n v_i (\xi_i) u_i \leq z)$ for a fixed stock allocation $u$ and fixed value of $z$.
Defining
\begin{equation}
\label{eq:s-F}
F(u, \xi) = -\sum_{i = 1}^n v_i(\xi_i) u_i,
\end{equation}
this is equivalent to estimating $\PP(F(u,\xi)\geq -z)$. %
Note that the function $F$ is concave in $\xi$ and linear in $u$, and in particular, it satisfies condition \ref{assume:concave_F} in \Cref{thm:ldt_verifiable}.

We downloaded 1,000 days' worth of daily NASDAQ stock price data for 50 stocks to estimate the drift parameter $\theta_i$, standard deviation $\sigma_i$ and the distribution for the random parameter $\xi$ (Gaussian or Gaussian mixture). Specifically, we first compute the compounded rates of $v_i$; \textit{i.e.}, the logarithm of the ratio of stock prices in two consecutive days. For a Gaussian distribution, we then compute the empirical mean and covariance corresponding to these rates to obtain $\theta_i-\sigma_i^2/2$ and covariance for $\xi_i$ respectively. For a Gaussian mixture, we use the \texttt{GaussianMixtures.jl} package to estimate the parameters of a $2$- and $3$-component Gaussian mixture using the expectation-maximization algorithm. %
The allocation $u$ is fixed and chosen randomly from a uniform distribution in $[0,1]^n$ subject to the constraint $\sum_{i=1}^n u_i =1$. %
We fix $T=10$, and compare the accuracy of estimating $\PP(F(u,\xi)\geq -z)$ for different fixed values of $z$ using: \textit{(i)} our first- ($P_1$) and second-order ($P_2$) LDT probability approximations, with \textit{(ii)} a Monte Carlo sample average approximation $P^{MC}_{N}(u,z)$ using $N$ independent samples.
We use $N=10^7$ samples to estimate the true probability as $P^{true} = P^{MC}_{10^7}(u,z)$, and define errors
\begin{equation*}
\label{eq:error}
\begin{aligned}
\epsilon^{LDT}_k(u,z)&:=\log_{10} P_k(u,z,\xistar) -\log_{10} P^{true}, \quad k \in \{1, 2\}, \\
\epsilon^{MC}_{N}(u,z)&:=\log_{10} P^{MC}_N(u,z) -\log_{10} P^{true}, \quad N \in \{10^4, 10^5\}.
 \end{aligned}
\end{equation*} 
We note that positive errors denote an over-estimation of the true probability, whereas negative values imply under-estimation.

\begin{table}[tbhp] 
    \footnotesize
        \caption{The comparison of the first- ($P_1$) and second-order ($P_2$) LDT probability estimates with sample average approximations $(P^{MC}_{N})$ for estimating $\PP(F(u,\xi)\geq -z)$ in the portfolio management problem. We drop the dependence on $u$ and $z$ for simplicity.
        }\label{tab:stock-G}
        \centering
				\begin{tabular}{cccccccccc}
				\toprule
				$z$	& $P^{true}$&$P_1$ & $P_2$&
		  $P^{MC}_{10^4}$ & $P^{MC}_{10^5}$ &
				$ \epsilon^{LDT}_1 $  & $\epsilon^{LDT}_2$ &
				 $ \epsilon^{MC}_{10^4}$ & $ \epsilon^{MC}_{10^5}$
				\\ 
				\midrule
                \multicolumn{10}{c}{\textit{Gaussian}} \\
                0.80 & $5.0\cdot 10^{-7}$ & $1.9\cdot 10^{-6}$ & $4.8\cdot 10^{-7}$ & $0.0$              & $0.0$              & $0.57$ & $-0.01$ & $-\infty$ & $-\infty$ \\
                0.82 & $5.4\cdot 10^{-6}$ & $2.2\cdot 10^{-5}$ & $6.6\cdot 10^{-6}$ & $0.0$              & $1.0\cdot 10^{-5}$ & $0.61$ & $ 0.08$ & $-\infty$ & $0.27$    \\
                0.84 & $5.8\cdot 10^{-5}$ & $1.9\cdot 10^{-4}$ & $6.3\cdot 10^{-5}$ & $0.0$              & $4.0\cdot 10^{-5}$ & $0.51$ & $ 0.04$ & $-\infty$ & $-0.16$   \\
                0.86 & $3.9\cdot 10^{-4}$ & $1.1\cdot 10^{-3}$ & $4.4\cdot 10^{-4}$ & $5.0\cdot 10^{-4}$ & $2.9\cdot 10^{-4}$ & $0.46$ & $ 0.05$ & $0.10$    & $-0.13$   \\
                \midrule
                \multicolumn{10}{c}{\textit{Gaussian mixture } $(M=2)$} \\
                0.80 & $6.6\cdot 10^{-5}$ & $2.2\cdot 10^{-4}$ & $7.5\cdot 10^{-5}$ & $0.0$              & $6.0\cdot 10^{-5}$ & $0.53$ & $0.06 $ & $-\infty$ & $-0.04$ \\
                0.82 & $1.8\cdot 10^{-4}$ & $5.7\cdot 10^{-4}$ & $2.1\cdot 10^{-4}$ & $3.0\cdot 10^{-4}$ & $2.0\cdot 10^{-4}$ & $0.49$ & $0.06 $ & $0.22$    & $0.04$ \\
                0.84 & $4.7\cdot 10^{-4}$ & $1.3\cdot 10^{-3}$ & $5.1\cdot 10^{-4}$ & $2.0\cdot 10^{-4}$ & $3.5\cdot 10^{-4}$ & $0.45$ & $0.04 $ & $-0.37$   & $-0.13$ \\
                0.86 & $1.1\cdot 10^{-3}$ & $3.0\cdot 10^{-3}$ & $1.1\cdot 10^{-3}$ & $1.2\cdot 10^{-3}$ & $1.3\cdot 10^{-3}$ & $0.42$ & $0.00 $ & $0.02$    & $0.05$ \\
                \midrule
                \multicolumn{10}{c}{\textit{Gaussian mixture } $(M=3)$} \\
                0.80 & $5.7\cdot 10^{-5}$ & $2.5\cdot 10^{-4}$ & $5.8\cdot 10^{-5}$ & $0.0$              & $4.0\cdot 10^{-5}$ & $0.65$ & $ 0.01$ & $-\infty$ & $-0.15$ \\
                0.82 & $2.1\cdot 10^{-4}$ & $7.9\cdot 10^{-4}$ & $2.0\cdot 10^{-4}$ & $1.0\cdot 10^{-4}$ & $2.2\cdot 10^{-4}$ & $0.57$ & $-0.02$ & $-0.32$   & $0.02$ \\
                0.84 & $6.5\cdot 10^{-4}$ & $2.1\cdot 10^{-3}$ & $6.1\cdot 10^{-4}$ & $1.0\cdot 10^{-3}$ & $7.9\cdot 10^{-4}$ & $0.51$ & $-0.03$ & $0.19$    & $0.09$ \\
                0.86 & $1.7\cdot 10^{-3}$ & $5.0\cdot 10^{-3}$ & $1.6\cdot 10^{-3}$ & $2.3\cdot 10^{-3}$ & $1.9\cdot 10^{-3}$ & $0.47$ & $-0.03$ & $0.13$    & $0.04$ \\
				\bottomrule
			\end{tabular} 
\end{table}

\Cref{tab:stock-G} shows that, as $z$ becomes smaller, the event $F(u,\xi)\geq -z$ becomes rarer and a sample average approximation with $N=10^4$ or $10^5$ samples is insufficient to accurately estimate the true probability.
For example we observe that for $z = 0.80$, the true probability is much smaller than $10^{-4}$;
in such cases, none of the $N=10^4$ samples record the rare event, on average, resulting in a probability estimate of $0$.
Stated differently, $P^{MC}_{10^4}$ and  $P^{MC}_{10^5}$ fail to approximate probabilities lower than $10^{-4}$ and $10^{-5}$, respectively, illustrating their high sample size requirement for accurate rare event estimation.

In contrast, our proposed LDT approximations are relatively less sensitive to the extremeness of the event (\textit{i.e.}, to values of $z$).
In particular, since $F$ is concave in $\xi$, \Cref{prop:concave_LDT_prob_upper_bound} ensures that the first-order approximation $P_1(u,z,\xistar)$ always overestimates the probability, and this is verified in \Cref{tab:stock-G}, since the errors $\epsilon^{LDT}_1$ are positive in all cases.
Nevertheless, the first-order estimates can be about three times the magnitude of the true probability, with log-errors of roughly $0.5$, across different $z$.
In contrast, the second-order estimates are highly accurate, with log-errors consistently less than $10^{-1}$, regardless of values of $z$, and they typically outperform the Monte Carlo estimates for both Gaussians and mixture distributions.
We highlight that, for rare events, it is much more important to correctly estimate the order of magnitude of the probability rather than its exact value; for example, it matters less that we estimate an event probability of $10^{-6}$ as $2\cdot 10^{-6}$, as compared to estimating it as $10^{-3}$.

\subsubsection{Optimal stock allocation}\label{sec:stock-prob}
We now turn to determining the optimal stock allocation.
We first define the value-at-risk (VaR) of a fixed allocation vector $u$:
\begin{equation}
\label{eq:s-VaR}
\VaR\nolimits_{1-\alpha}(u):=\argmax\limits_z\left\lbrace \PP\left(\sum\limits_{i = 1}^n v_i(\xi_i) u_i \geq z\right)\geq 1 - \alpha \right\rbrace.
\end{equation}
Our goal is to maximize the VaR of the total return at level $1-\alpha$, and consider $\alpha\ll 1$ reflecting a conservative investment approach.
This problem can be reformulated as the following variant of the chance-constrained optimization problem \eqref{eq:ccp}  using the definition of $F$ in \eqref{eq:s-F}:
\begin{equation}
\label{eq:portfolio_ccp}
\begin{aligned}
& \underset{u \in \mathbb{R}^n_{+}, z \in \RR}{\text{maximize}}
& & z \\
& \text{subject to}
& &  \quad \sum_{i=1}^n u_i = 1 \\
&&& \quad \PP\left(F(u,\xi) \geq -z\right)\leq \alpha. %
\end{aligned}
\end{equation}

For $T = 10$ and $\alpha = 10^{-4}$, \Cref{tab:stockT10a4-G} compares the computational performance of our LDT chance-constrained formulations with the smooth sigmoidal sample average formulation \eqref{eq:op-SA2} (with $\nu = \tau = 1$) and the CVaR sample average formulation \eqref{eq:op-CVaR}, each with $N = 10^5$ samples.
Here, $u^\star$ denotes the optimal allocation determined by each method, $z^\star$ denotes the corresponding objective value of \eqref{eq:portfolio_ccp}, and $\VaR_{1-\alpha}(u^\star)$ denotes the Value-at-Risk \eqref{eq:s-VaR} of the corresponding allocation $u^\star$ which we would like to be maximized.
As before, we estimate the true probability $P^{true}(u^\star,z^\star) = P^{MC}_{10^7}(u^\star,z^\star)$ using $N = 10^7$ Monte Carlo samples.
Therefore, the column $\log_{10} P^{true} (u^\star,z^\star)-\log_{10}  \alpha$ shows the extent of feasibility of the chance constraint in each case, where negative values indicate satisfaction of the chance constraint whereas positive values indicate that the corresponding $u^\star$ and $z^\star$ fail to satisfy the chance constraint in \eqref{eq:portfolio_ccp}.
We also report solution times of the different methods. 

\begin{table}[tbhp] 
	{\footnotesize
		\caption{Comparison of solution methods for the optimal stock allocation problem \eqref{eq:portfolio_ccp} with $T=10$ and $\alpha=10^{-4}$.
		}\label{tab:stockT10a4-G}
		\begin{center}
			\begin{tabular}{cccccc}
				\toprule 
				Method	& $z^\star$ & $\VaR_{1-\alpha}(u^\star)$&$P^{true}(u^\star,z^\star)$& $\log_{10}  P^{true}(u^\star,z^\star)-\log_{10}  \alpha$& Time [sec] \\ \midrule
                \multicolumn{6}{c}{\textit{Gaussian}} \\
				LDT $(k = 1)$     & $9.53\cdot 10^{-1}$ & $9.54\cdot 10^{-1}$ & $7.7\cdot 10^{-5}$ & $-0.11  $ & $0.5  $ \\ 
				LDT $(k = 2)$     & $9.54\cdot 10^{-1}$ & $9.54\cdot 10^{-1}$ & $9.5\cdot 10^{-5}$ & $-0.02  $ & $197.0$ \\ 
				SAA $(N = 10^5)$  & $9.02\cdot 10^{-1}$ & $9.53\cdot 10^{-1}$ & $0.0$              & $-\infty$ & $638.7$ \\ 
				CVaR $(N = 10^5)$ & $9.53\cdot 10^{-1}$ & $9.52\cdot 10^{-1}$ & $1.5\cdot 10^{-4}$ & $0.17   $ & $527.6$ \\ \midrule
                \multicolumn{6}{c}{\textit{Gaussian mixture } ($M = 2$)} \\
                LDT $(k = 1)$     & $9.15\cdot 10^{-1}$ & $9.20\cdot 10^{-1}$ & $5.4\cdot 10^{-5}$ & $-0.26$ & $11.3$ \\ 
                LDT $(k = 2)$     & $9.21\cdot 10^{-1}$ & $9.23\cdot 10^{-1}$ & $8.3\cdot 10^{-5}$ & $-0.08$ & $85.7$ \\ 
                SAA $(N = 10^5)$  & $9.20\cdot 10^{-1}$ & $9.18\cdot 10^{-1}$ & $1.2\cdot 10^{-4}$ & $0.09$  & $1,147.4$ \\ 
                CVaR $(N = 10^5)$ & $9.19\cdot 10^{-1}$ & $9.17\cdot 10^{-1}$ & $1.2\cdot 10^{-4}$ & $0.07$  & $420.4$ \\
				\bottomrule 
			\end{tabular} 
		\end{center}
	}
\end{table}

From \Cref{tab:stockT10a4-G}, we observe that our proposed second-order LDT formulation consistently attains the maximum VaR.
Moreover, its optimal objective value $z^\star$ is quite close to the true VaR, indicating that the second-order probability estimate is also accurate.
Notably, these improved objective values also come at a much lower solution time compared to the sampling-based methods.
The first-order LDT formulation attains a slightly smaller VaR compared to the second-order formulation, but it is still better than the sampling-based methods.
Moreover, the solutions $(u^\star, z^\star)$ obtained by the sampling-based methods are not feasible in the original problem \cref{eq:portfolio_ccp}, since the corresponding $\log_{10}  P^{true}(u^\star,z^\star)-\log_{10}  \alpha > 0$, implying that we need to increase the number of samples $N$ for these methods to obtain a feasible solution.
However, given that their solution times are already quite excessive with $N=10^5$ samples; further increasing $N$ will only amplify these times even more.
In contrast, the first- and second-order LDT formulations solve in $\sim$10 and $\sim$100 seconds respectively, with their solution times being largely unaffected by the choice of $\alpha$.
In conclusion, the LDT formulations outperform the sampling-based methods both in terms of objective value and solution time.
The choice between the first- and second-order formulations can be made on a case-by-case basis; whereas the former is more computationally efficient, the latter is more accurate, in general.

\begin{figure}[tbhp]\centering 
	\begin{tikzpicture}
	\begin{axis}[width=0.9\textwidth, height=0.4\textwidth, 
	scale only axis, ymax=0.91, ymin=-1e-8,xmin=1, xmax=50,
	legend style = {font=\small,nodes=right},
	xlabel=Symbol of stock,
	ylabel=Optimal allocation $u^\star$,
	x tick label style={font=\tiny},
	ybar interval =1,
	legend pos = north west,
	x label style={at={(axis description cs:0.5,-0.1)},anchor=north},
	xtick={1,2,...,50},
	xticklabels={AAL ,AAME ,AAON ,AAWW ,AAXJ ,AAXN ,ABCB ,ABIO ,ABMD ,ABTX ,ABUS ,ACBI ,ACGL ,ACHC ,ACIA ,ACIW ,ACNB ,ACRX ,ACTG ,ACWI ,ACWX ,ADAP ,ADES ,ADI ,ADMA ,ADP ,ADSK ,ADUS ,ADVM ,AEGN ,AEHR ,AEIS ,AEMD ,AERI ,AEY ,AEYE ,AFH ,AFMD ,AGEN ,AGFS ,AGFSW ,AGIO ,AGLE ,AGRX ,AGTC ,AGYS ,AGZD ,AHPI ,AIA ,AIMC 
	},
	x tick label style={rotate=90,anchor=east},
	]
	
	\addplot +[draw=\cLDTF!50!black,fill=\cLDTF!50]
	coordinates {
		(1, 2e-08) (2, 5e-07) (3, 7e-09) (4, 5e-06) (5, 4e-07)  
		(6, 3e-08) (7, 3e-08) (8, 5e-07) (9, 3e-08) (10, 0.04)  
		(11, 2e-08) (12, 0.002) (13, 2e-07) (14, 4e-09) (15, 0.04)  
		(16, 4e-08) (17, 0.06) (18, 1e-08) (19, 0.02) (20, 6e-08)  
		(21, 1e-07) (22, 0.005) (23, 4e-09) (24, 3e-06) (25, 5e-07)  
		(26, 6e-06) (27, 2e-08) (28, 0.03) (29, 2e-08) (30, 6e-08)  
		(31, 4e-08) (32, 0.02) (33, 3e-08) (34, 6e-07) (35, 0.003)  
		(36, 0.01) (37, 0.0007) (38, 6e-07) (39, 3e-08) (40, 3e-09)  
		(41, 0.002) (42, 2e-08) (43, 3e-06) (44, 3e-08) (45, 5e-07)  
		(46, 1e-08) (47, 0.8) (48, 0.008) (49, 1e-07) (50, 5e-09)
	};
	\addlegendentry{LDT $(k = 1)$ }
	\addplot +[draw=\cLDTS!50!black,fill=\cLDTS!50]
	coordinates {
	(1, -1e-08) (2, -9e-09) (3, -1e-08) (4, 0.01) (5, -9e-09)  
	(6, -1e-08) (7, -1e-08) (8, 0.003) (9, -1e-08) (10, 0.04)  
	(11, -1e-08) (12, 0.02) (13, -9e-09) (14, -1e-08) (15, 0.05)  
	(16, -1e-08) (17, 0.05) (18, -1e-08) (19, 0.02) (20, -1e-08)  
	(21, -1e-08) (22, 0.007) (23, -1e-08) (24, -8e-09) (25, -8e-09)  
	(26, 3e-08) (27, -1e-08) (28, 0.04) (29, -1e-08) (30, -1e-08)  
	(31, -1e-08) (32, -4e-10) (33, -1e-08) (34, -8e-09) (35, 0.007)  
	(36, 0.006) (37, 0.005) (38, -9e-09) (39, -1e-08) (40, -1e-08)  
	(41, 0.006) (42, -1e-08) (43, -9e-09) (44, -1e-08) (45, 0.005)  
	(46, -1e-08) (47, 0.7) (48, 0.01) (49, -1e-08) (50, -1e-08)  
	};
	\addlegendentry{LDT $(k = 2)$ }

		\addplot+[draw=\cSATfive!50!black, fill=\cSATfive!70!black] coordinates {
			(1, 0.006) (2, 0.009) (3, 7e-09) (4, 0.02) (5, 4e-07)  
			(6, 1e-07) (7, 7e-08) (8, 0.0007) (9, 2e-08) (10, 0.04)  
			(11, 2e-08) (12, 2e-07) (13, 1e-07) (14, 6e-09) (15, 0.03)  
			(16, 4e-08) (17, 0.06) (18, 2e-08) (19, 0.008) (20, 6e-08)  
			(21, 1e-07) (22, 0.004) (23, 3e-09) (24, 2e-07) (25, 8e-08)  
			(26, 1e-06) (27, 2e-08) (28, 1e-07) (29, 3e-08) (30, 2e-08)  
			(31, 2e-08) (32, 7e-07) (33, 1e-08) (34, 2e-07) (35, 0.04)  
			(36, 2e-07) (37, 0.005) (38, 9e-08) (39, 0.005) (40, 6e-09)  
			(41, 0.005) (42, 2e-08) (43, 7e-08) (44, 1e-06) (45, 0.02)  
			(46, 1e-08) (47, 0.7) (48, 0.02) (49, 2e-07) (50, 4e-09) 
	};
		\addlegendentry{SAA $(N = 10^5) $ }

	\addplot+[draw=\cCVaRfive!50!black, fill=\cCVaRfive!50] 
	coordinates {
		(1, 3e-08) (2, 0.003) (3, 6e-09) (4, 0.02) (5, 2e-07)  
		(6, 7e-08) (7, 3e-08) (8, 0.002) (9, 2e-08) (10, 0.07)  
		(11, 1e-08) (12, 2e-07) (13, 1e-07) (14, 5e-09) (15, 0.05)  
		(16, 4e-08) (17, 0.06) (18, 2e-08) (19, 2e-07) (20, 5e-08)  
		(21, 7e-08) (22, 0.007) (23, 4e-09) (24, 3e-06) (25, 3e-08)  
		(26, 3e-07) (27, 2e-08) (28, 0.01) (29, 1e-08) (30, 4e-08)  
		(31, 2e-08) (32, 2e-06) (33, 2e-08) (34, 6e-07) (35, 0.04)  
		(36, 7e-08) (37, 2e-07) (38, 6e-08) (39, 0.001) (40, 6e-09)  
		(41, 0.0007) (42, 3e-08) (43, 2e-07) (44, 4e-08) (45, 0.006)  
		(46, 1e-08) (47, 0.7) (48, 0.02) (49, 8e-08) (50, 5e-09) 
	};
	\addlegendentry{CVaR $(N = 10^5) $ }
	\end{axis}
	\end{tikzpicture}
	\caption{Optimal stock allocation determined by different methods with $T=10$ and $\alpha=10^{-4}$, for the case where $\xi$ follows a Gaussian mixture distribution with two components.}\label{fig:stockT10a4-GM}
\end{figure}
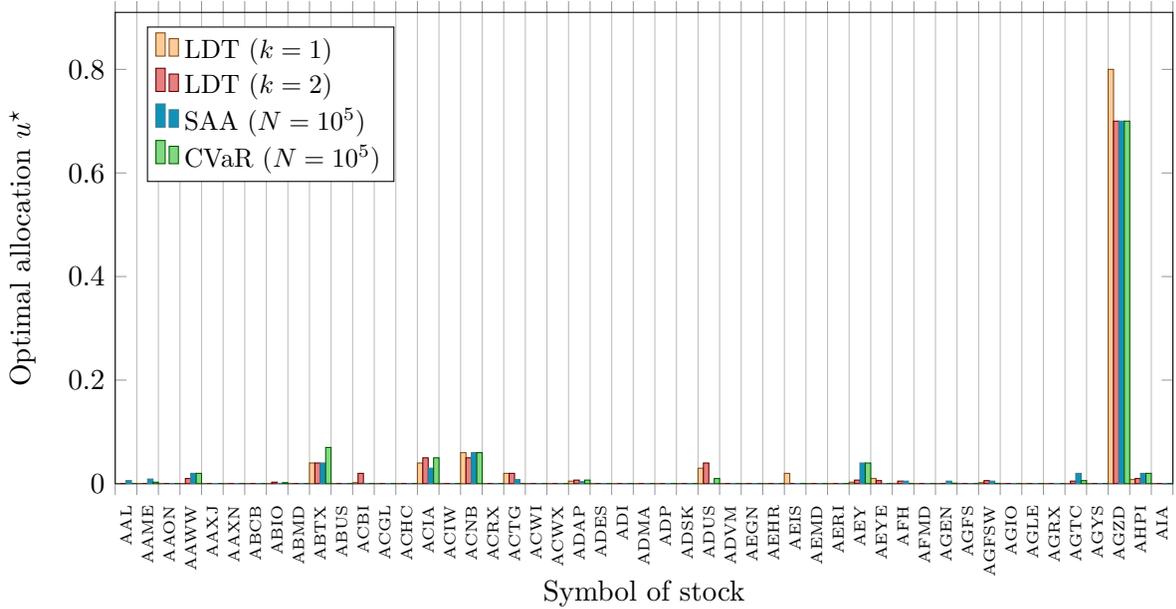

\Cref{fig:stockT10a4-GM} shows the optimal stocks allocations obtained by different methods in the case of the Gaussian mixture.
Interestingly, all methods indicate that we should invest majority of the initial wealth in the stock ``AGZD'',and never invest in some of the other stocks.
However, the first-order LDT formulation determines a slightly different allocation for this stock compared with other methods, implying that we should use the second-order LDT allocation to obtain higher returns.
Curiously, a closer look at \Cref{tab:stockT10a4-G} reveals that the VaR of the total return is less than the initial wealth ($=1$) in all cases, indicating that it may not be wise to invest at all over the short time period $T=10$ with an overly conservative risk threshold $\alpha = 10^{-4}$. %

\subsection{Short column design}\label{sec:shortcolumn} 
In the last application, the function $F$ was linear in the decisions $u$. 
We now consider a nonlinear model for the design of a short column \cite{chaudhuri2020risk,aoues2010benchmark} with uncertain material properties, and subject to a chance constraint on material failure. %

Consider a short column with a rectangular cross-section of width $w$ and height $h$, %
that is subject to an uncertain axial load $\xi_F$, bending moment $\xi_M$ and an uncertain yield stress $\xi_Y$.
We model these as a three-dimensional random vector $\xi = [\xi_F, \xi_M, \xi_Y]$.
Our goal is to determine $w$ and $h$ so that the area of the cross-section $wh$ is minimized, and the design should not experience material failure with very high probability of $1-\alpha$, where $\alpha \in (0,1) $ is some given threshold.
The parameter-to-observation map $F$ can be derived from the elastic-plastic constitutive law, and it has a general nonlinear dependence on both the decisions $u = (w, h) \in \RR^2$ as well as the uncertainty $\xi \in \RR^3$:
\begin{equation}
\label{eq:sc-F}
F(u,\xi):=\frac{4\xi_M}{wh^2\exp(\xi_Y)}+\frac{\xi_F^2}{w^2h^2\exp(2\xi_Y)}.
\end{equation}
The optimization problem is
\begin{equation}
\label{eq:sc-op}
\begin{aligned}
& \underset{u=(w,h) \in\RR^2}{\text{minimize}}
& & wh \\
& \text{subject to}
& & L_w\leq w\leq U_w, \quad L_h\leq h\leq U_h, \\
& & & \PP(F(u,\xi)\geq 1)\leq \alpha, %
\end{aligned}
\end{equation}
where $L_w,L_h$ and $U_w,U_h$ are given lower and upper bounds on $w$ and $h$, respectively, and we set $[L_w, U_w, L_h, U_h] = [5,15, 15,25]$.

\subsubsection{Gaussian random parameter} We first consider the random vector $\xi\sim\Gauss{\mu}{\Sigma}$ with mean and covariance given as follows, adapted from the example in \cite{chaudhuri2020risk}:
\begin{equation}
\label{eq:sc-dis}
\mu=\left[500, 2000, 1.604\right],\quad \Sigma = \left[ \begin{matrix}
10000 & 20000 &0\\
20000 & 160000 & 0\\
0 &0&  0.00995
\end{matrix}\right] .
\end{equation}
The covariance matrix indicates that the correlation coefficient between the axial force $\xi_F$ and bending moment $\xi_M$ is 0.5, whereas the yield stress $\exp(\xi_Y)$ is uncorrelated to them and follows a log-normal distribution with mean $5$ and standard deviation $0.5$.

\begin{figure}[tbhp]\centering
	\begin{tikzpicture}[]
	\begin{groupplot}[group style = {columns=3, 
	    rows=1,
		vertical sep =70pt, horizontal sep = 30pt}, width = 4.5cm, height = 5.0cm, title style = {font= \small},label style={font=\small},tick label style={font=\footnotesize}, ,xmax=2*1e-1,xmin=0.5*1e-6, cycle list name= my color]
	\nextgroupplot[xmode=log,
	scale only axis,
	xlabel={Risk threshold $\alpha$},
	title = {Optimal area $wh$ [$\text{m}^2$]},
	legend style = {font=\small,nodes=right, legend to name=grouplegend},
	legend columns=4]
	
	\addplot table[x=alpha,y=obj] {\data/sc_FO.txt};
	\addlegendentry{LDT $(k = 1)$}
	
	\addplot table[x=alpha,y=obj] {\data/sc_SO.txt};
	\addlegendentry{LDT $(k = 2)$}
	
	\addplot table[x=alpha,y=obj] {\data/sc_CVaR.txt};
	\addlegendentry{CVaR $(N =10^5)$}
	\addplot table[x=alpha,y=obj] {\data/sc_SA2_K1m200logN2.txt};
	\addlegendentry{SAA $(N =10^2)$}
	
	\addplot table[x=alpha,y=obj] {\data/sc_SA2_K1m200logN3.txt};
	\addlegendentry{SAA $(N =10^3)$}
	
	\addplot table[x=alpha,y=obj] {\data/sc_SA2_K1m200logN4.txt};
	\addlegendentry{SAA $(N =10^4)$}
	
	\addplot table[x=alpha,y=obj] {\data/sc_SA2_K1m200logN5.txt};
	\addlegendentry{SAA $(N =10^5)$}

	\nextgroupplot[xmode=log, scale only axis,
	xlabel={Risk threshold $\alpha$},
	title={ $\log_{10}\PP(F(u^\star,\xi)\geq 1)-\log_{10} \alpha$},ymax=4.5,ymin=-1]
	
	\addplot table[x=alpha,y expr=\thisrowno{6}/2.3] {\data/sc_FO.txt};
	\addplot table[x=alpha,y expr=\thisrowno{6}/2.3] {\data/sc_SO.txt};
	\addplot table[x=alpha,y expr=\thisrowno{6}/2.3] {\data/sc_CVaR.txt};
	\addplot table[x=alpha,y expr=\thisrowno{6}/2.3] {\data/sc_SA2_K1m200logN2.txt};
	\addplot table[x=alpha,y expr=\thisrowno{6}/2.3] {\data/sc_SA2_K1m200logN3.txt};
	\addplot table[x=alpha,y expr=\thisrowno{6}/2.3] {\data/sc_SA2_K1m200logN4.txt};
	\addplot table[x=alpha,y expr=\thisrowno{6}/2.3] {\data/sc_SA2_K1m200logN5.txt};
	
	\addplot[name path=zero, color=black!50,line width=0.1pt] coordinates {(0.5*1e-7,0) (2*1e-1,0)};
	\addplot[name path=under, color=black!50,line width=0.1pt, densely dotted] coordinates {(0.5*1e-6,-1) (2*1e-1,-1)};
	\addplot [
	thick,
	color=blue,
	fill=cyan, 
	fill opacity=0.05
	]
	fill between[
	of=zero and under,
	soft clip={domain=0.5*1e-6:2*1e-1},
	];

	\nextgroupplot[xmode=log, ymode=log, scale only axis,
	xlabel={Risk threshold $\alpha$},
	title={Time [sec]}]
	\addplot table[x=alpha,y=t] {\data/sc_FO.txt};
	\addplot table[x=alpha,y=t] {\data/sc_SO.txt};
	\addplot table[x=alpha,y=t] {\data/sc_CVaR.txt};
	\addplot table[x=alpha,y=t] {\data/sc_SA2_K1m200logN2.txt};
	\addplot table[x=alpha,y=t] {\data/sc_SA2_K1m200logN3.txt};
	\addplot table[x=alpha,y=t] {\data/sc_SA2_K1m200logN4.txt};
	\addplot table[x=alpha,y=t] {\data/sc_SA2_K1m200logN5.txt};
	\end{groupplot}
	\node[black] at ($(group c2r1) + (-0.0cm,4.2cm)$) {\pgfplotslegendfromname{grouplegend}}; 
	\end{tikzpicture}
	\caption{Comparison of different methods for solving the short column design problem with Gaussian random parameter. The parameters for SAA are $\nu=200$ and $\tau=200$. The left figure shows the optimal area (objective value) for each methods. The middle one shows the feasibility of the solutions, we expect the values $\log_{10} \PP(F(u^\star,\xi)\geq 1)-\log_{10} \alpha$ to be negative (fall in the cyan area), meaning the solution $u^\star$ satisfies the chance constraint. The right figure shows the time costs for different methods.
	}\label{fig:sc-G}
\end{figure}
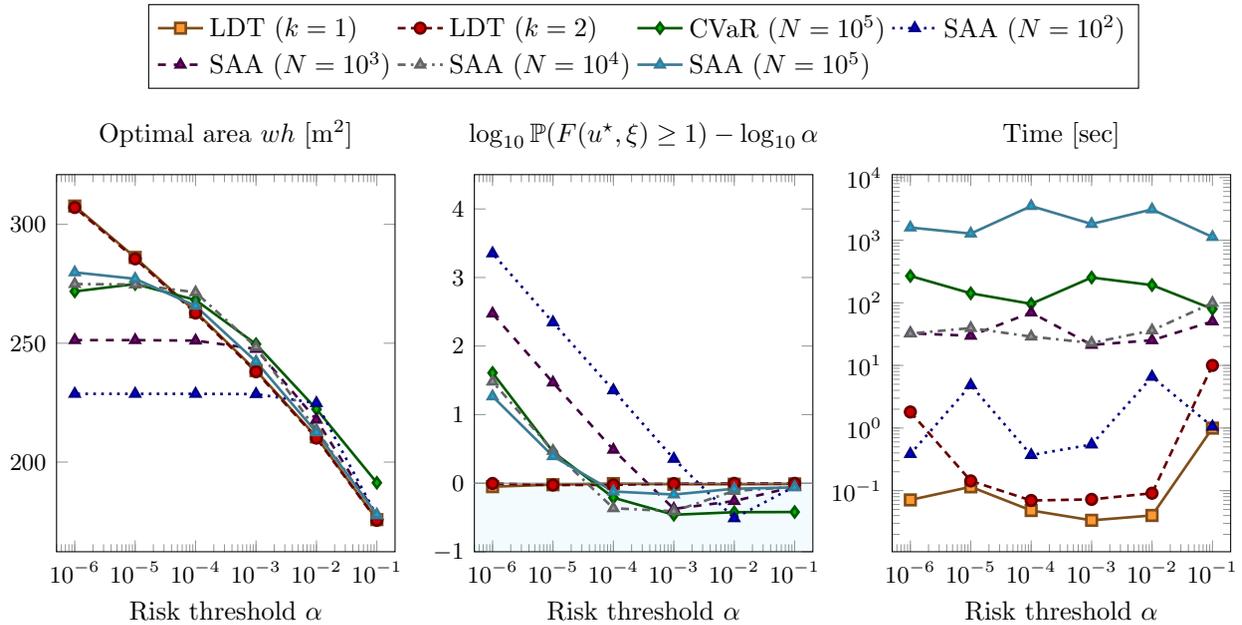

In \Cref{fig:sc-G}, we compare the performance of different methods for solving the short column design problem with Gaussian random parameter \eqref{eq:sc-dis}, as a function of the risk threshold $\alpha$ varying from $10^{-1}$ to $10^{-6}$. %
The left-most figure shows that our LDT chance constrained formulations provide similar optimal solutions, indicating that $F$ is near-linear.
Moreover, when $\alpha$ becomes smaller, although the optimal areas determined by the sampling-based methods appear to be smaller than the ones determined by the LDT approaches (near-horizontal parts in the left-most figure),
the middle figure shows that these solutions are, in fact, not feasible because they fail to satisfy the probability constraint (since their curves fall outside the negative part highlighted in blue).
In fact, as far as the sampling-based methods go, the solutions determined using $N$ samples can give reliable solutions only for $\alpha\geq 10^{-(N-1)}$. %
In contrast, the middle figure shows that the solutions determined by the LDT formulations always satisfy the chance constraint for all values of $\alpha$, and their corresponding optimal areas coincide with those of the sampling-based methods, whenever the latter are also feasible.
In fact, they are able to find feasible designs even for extreme cases where sampling-based methods with $N\leq 10^5$ samples fail.
Finally, the right-most figure shows that the computational time increases severely for sampling-based methods as the number of samples $N$ increases, since the problem size in \eqref{eq:op-SA2} and \eqref{eq:op-CVaR} scales as $O(N)$.
The middle figure suggests that, for very high reliability $\alpha=10^{-6}$, we need to increase $N$ to obtain a feasible solution, and thus the corresponding time is also expected to increase.
The LDT approaches are much more efficient, with both the first- and second-order formulations solving in roughly $10\%$ of the time compared to SAA with $N=10^2$ samples. %

\subsubsection{Gaussian mixture random parameter} 
\begin{figure}[tbhp]\centering
	\begin{tikzpicture}[]
	\begin{groupplot}[group style = {columns=3, 
	    rows=1,
		vertical sep =70pt, horizontal sep = 30pt}, width = 4.5cm, height = 5.0cm, title style = {font= \small},label style={font=\small},tick label style={font=\footnotesize}, ,xmax=2*1e-1,xmin=0.5*1e-6, cycle list name= my color]
	\nextgroupplot[xmode=log, 
	scale only axis,
	xlabel={Risk threshold $\alpha$},
	title = {Optimal area $wh$ [$\text{m}^2$]},
	legend style = {font=\small,nodes=right, legend to name=grouplegend2},legend columns=4]
	
	\addplot table[x=alpha,y=obj] {\data/sc_FO_K2.txt};
	\addlegendentry{LDT $(k = 1)$}
	\addplot table[x=alpha,y=obj] {\data/sc_SO_K2.txt};
	\addlegendentry{LDT $(k = 2)$}

	\addplot table[x=alpha,y=obj] {\data/sc_CVaR_K2.txt};
	\addlegendentry{CVaR $(N =10^5)$}
	
	\addplot table[x=alpha,y=obj] {\data/sc_SA2_K2m200logN2.txt};
	\addlegendentry{SAA $(N =10^2)$}
	
	\addplot table[x=alpha,y=obj] {\data/sc_SA2_K2m200logN3.txt};
	\addlegendentry{SAA $(N =10^3)$}
	
	\addplot table[x=alpha,y=obj] {\data/sc_SA2_K2m200logN4.txt};
	\addlegendentry{SAA $(N =10^4)$}
	
	\addplot table[x=alpha,y=obj] {\data/sc_SA2_K2m200.txt};
	\addlegendentry{SAA $(N =10^5)$}

	\nextgroupplot[xmode=log, scale only axis,
	xlabel={Risk threshold $\alpha$},
	title={$\log_{10}\PP(F(u^\star,\xi)\geq 1)-\log_{10} \alpha$},ymin=-1,ymax=4.5]
	\addplot table[x=alpha,y expr=\thisrowno{6}/2.3] {\data/sc_FO_K2.txt};
	\addplot table[x=alpha,y expr=\thisrowno{6}/2.3] {\data/sc_SO_K2.txt};
	\addplot table[x=alpha,y expr=\thisrowno{6}/2.3] {\data/sc_CVaR_K2.txt};
	\addplot table[x=alpha,y expr=\thisrowno{6}/2.3] {\data/sc_SA2_K2m200logN2.txt};
	\addplot table[x=alpha,y expr=\thisrowno{6}/2.3] {\data/sc_SA2_K2m200logN3.txt};
	\addplot table[x=alpha,y expr=\thisrowno{6}/2.3] {\data/sc_SA2_K2m200logN4.txt};
	\addplot table[x=alpha,y expr=\thisrowno{6}/2.3] {\data/sc_SA2_K2m200.txt};
	
	\addplot[name path=zero, color=black!50,line width=0.1pt] coordinates {(0.5*1e-7,0) (2*1e-1,0)};
	\addplot[name path=under, color=black!50,line width=0.1pt, densely dotted] coordinates {(0.5*1e-6,-1) (2*1e-1,-1)};
	\addplot [
	thick,
	color=blue,
	fill=cyan, 
	fill opacity=0.05
	]
	fill between[
	of=zero and under,
	soft clip={domain=0.5*1e-6:2*1e-1},	];

	\nextgroupplot[xmode=log, ymode=log, scale only axis,
	xlabel={Risk threshold $\alpha$},
	title={Time [sec]}]
	\addplot table[x=alpha,y=t] {\data/sc_FO_K2.txt};
	\addplot table[x=alpha,y=t] {\data/sc_SO_K2.txt};
	\addplot table[x=alpha,y=t] {\data/sc_CVaR_K2.txt};
	\addplot table[x=alpha,y=t] {\data/sc_SA2_K2m200logN2.txt};
	\addplot table[x=alpha,y=t] {\data/sc_SA2_K2m200logN3.txt};
	\addplot table[x=alpha,y=t] {\data/sc_SA2_K2m200logN4.txt};
	\addplot table[x=alpha,y=t] {\data/sc_SA2_K2m200.txt};
	\end{groupplot}
	\node[black] at ($(group c2r1) + (-0.0cm,4.2cm)$) {\pgfplotslegendfromname{grouplegend2}}; 
	\end{tikzpicture}
	\caption{Comparison of different methods for solving the short column design problem with Gaussian mixture random parameter with two components. The parameters for SAA are $\nu=200$ and $\tau=200$.}\label{fig:sc-GM}
\end{figure}
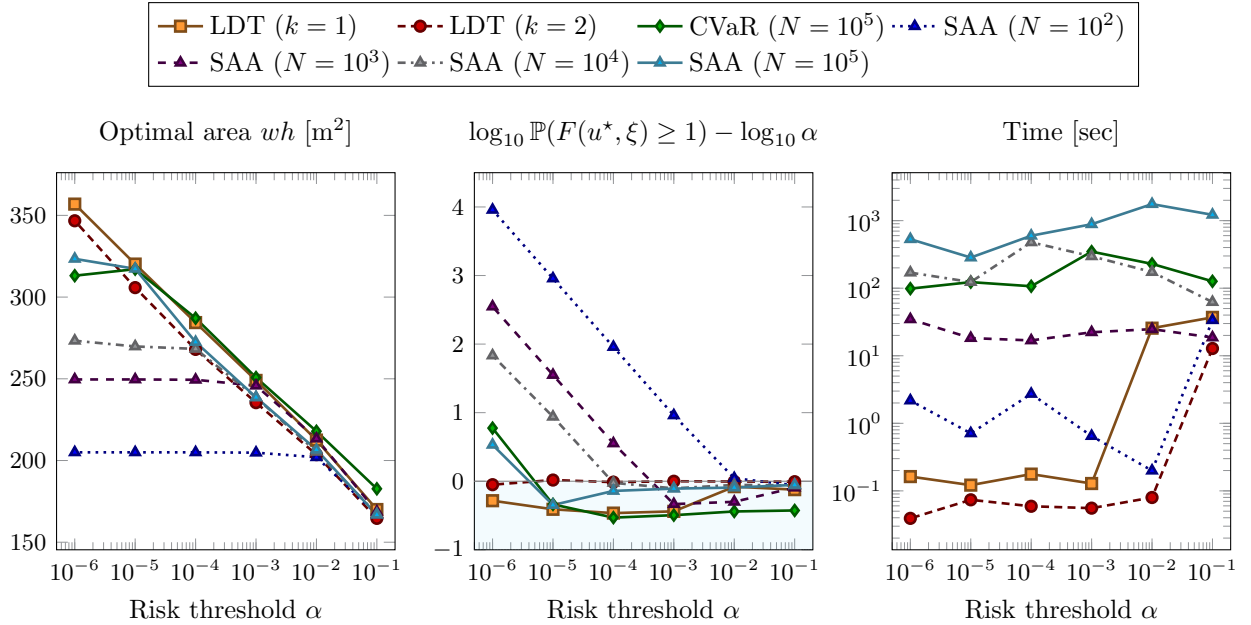

We now consider $\xi\sim 0.5\Gauss{\mu_1}{\Sigma_1}+ 0.5\Gauss{\mu_2}{\Sigma_2}$ to follow a Gaussian mixture distribution with two components, where the mean $\mu_1$ and covariance $\Sigma_1$ of the first Gaussian component have the same values as in \eqref{eq:sc-dis}, whereas 
the mean and covariance of the second Gaussian component is given as follows:
\begin{equation}
\mu_2=\left[100, 1000,1.0849\right], \quad
\Sigma_2 = \left[ 
\begin{matrix}
10000& 20000 &0\\
20000 &160000 &0\\
0 &0 &0.0274
\end{matrix}\right] .
\end{equation}
The results are summarized in \Cref{fig:sc-GM}, and are largely similar to the case of Gaussians. The main differences are the different optimal areas determined by the two LDT approaches, with the second-order formulation yielding a slightly better objective value while still obtaining a feasible design.
As before, both LDT approaches are significantly more efficient compared with sampling-based methods in terms of computational time and sensitivity to event extremeness.

\subsection{Optimal boundary control for steady-state advection-diffusion}\label{sec:PDE}
In the previous sections, we showed the performance of our proposed formulations approach for problems that are nonlinear in the uncertainty, but still admit closed-form expressions of the parameter-to-observation map $F$. In this section, we consider a much more complicated PDE system, where we do not have closed-form expressions.
In particular, these represent examples where it is infeasible to use sampling-based methods, since they would need to incorporate as many PDE solves as the number of samples.
In contrast, we demonstrate that our LDT approaches continue to provide satisfactory results.

Consider a 2D steady-state advection-diffusion equation in the unit square $\Omega=[0,1]\times[0,1]$, and denote the boundary $\Gamma_c={0}\times[0,1]$ and $\Gamma_n=\partial\Omega\backslash\Gamma_c$. 
\begin{subequations}\label{eq:sad-PDE}
	\begin{alignat}{3}
	-\nabla\cdot \left( \kappa(x,\xi)\nabla y(x)\right) + w(x)\cdot \nabla y(x) & = f(x,\xi), & x\in\Omega,\\
	\left(\kappa(x,\xi)\nabla y(x)\right) \cdot n(x) &= \frac{1}{\epsilon_0} \left(u(x)-y(x)\right), & \text{on }\Gamma_c, \label{eq:Dbd}\\
	\left(\kappa(x,\xi)\nabla y(x)\right) \cdot n(x) & = 0, & \text{on }\Gamma_n.
	\end{alignat}
\end{subequations}
Here, $y$ is the state representing the temperature at each location $x$ and $u$ is the Dirichlet boundary control on the boundary $\Gamma_c$.
We enforce $y(x)=u(x)$ on $\Gamma_c$ through a Robin-type condition \eqref{eq:Dbd} to avoid regularity limitations of Dirichlet boundary control.
Also, $\epsilon_0=10^{-4}$ is a penalty parameter, $w=[1,0]^\top$ is a fixed velocity field, and $n(x)$ represents the outward-pointing normal at the boundary.
The diffusion coefficient $\kappa(x,\xi)$ and force $f(x,\xi)$ are subject to random uncertainty $\xi$ with the following nonlinear dependence:
\begin{equation}
\label{eq:sad-kappa}
\begin{aligned}
\kappa(x,\xi) &= \left\lbrace \begin{array}{ll}
0.8,& x\in [0,1]\times[0.6,1] , \\ 
\exp(\xi_1),&x\in [0,1]\times[0,0.6]
\end{array} \right. \\
f(x,\xi) & =20\exp\left(-\frac{(x_1-\xi_2)^2}{0.1}\right)\exp\left(-\frac{(x_2-0.5)^2}{0.1}\right).
\end{aligned}
\end{equation}

Observe that, since the temperature $y$ is a solution of \eqref{eq:sad-PDE}, it is uncertain due to the randomness of $\xi$.
In applications like hyperthermia treatment for cancer therapy \cite{deuflhard2011mathematical}, it is important to ensure that the (uncertain) temperature in some areas (\textit{e.g.}, essential organs) do not exceed some critical threshold to protect organs from thermal damage.
Motivated by these applications, our goal is to ensure that the average temperature in the sub-domain $\Omega_0=[0.4,0.6]\times[0.4,0.6]$ remains less than $z$ with high probability $1-\alpha$.
This leads to the following definition of the parameter-to-observation $F(u,\xi)$ as the average temperature in $\Omega_0$ for a given realization of the random parameter $\xi$ and control vector $u$:
\begin{equation}
\label{eq:sad-F}
F(u,\xi) := \frac{1}{|\Omega_0|} \int_{\Omega_0}y(x;u,\xi) \mathrm{d}x,
\end{equation}
where $y(x;u,\xi)$ is the solution of PDE \eqref{eq:sad-PDE} at location $x$ given $\xi$ and $u$.
Our goal is to control the temperature at the boundary $\Gamma_c$ as close to 0 as possible, while ensuring that the average temperature in $\Omega_0$ does not exceed $z$ with high probability:
\begin{equation}
\label{eq:sad-minreg}
\begin{aligned}
& \underset{u}{\text{minimize}}
& & \frac{1}{2} \int_{\Gamma_c}u^2(x) dx,\\
& \text{subject to}
 & & \PP(F(u,\xi)\geq z)\leq \alpha,
\end{aligned}
\end{equation}
where $z$ is a given upper bound on the average temperature in $\Omega_0$, and $\alpha\ll 1$.

Unlike the other applications, note that we do not have a closed-form expression for $F( u,\xi)$ with respect to $u$ or $\xi$, and we can only access it through the solution of the PDE in \eqref{eq:sad-PDE}.
Nevertheless, it can still be characterized implicitly as a linear system upon discretization of the PDE.
If we denote the dimension of the design $u$ to be $m$ and use a finite difference method to dicretize the PDE, then the resulting linear system has size $m^2\times m^2$.
Since the complexity of solving one linear problem is $O(m^2)$, 
resorting to a sampling-based methods will require us to model $N$ realizations of $F$ as constraints in the optimization problem, each involving a PDE-based linear solve.
This would result in a very large-scale problem even for modest numbers of samples and control inputs, making sampling-based methods practically inefficient.
This motivates alternate LDT approaches.

\begin{table}[tbhp] 
	{\footnotesize
		\caption{Optimal objective values and constraint feasibilities of the boundary control problem \eqref{eq:sad-minreg} with $z=1$ for different $\alpha$'s using our first- and second order LDT chance-constrained formulations.} \label{tab:sad1}
		\begin{center}
			\begin{tabular}{cccc}
				\toprule 
				$\alpha$ &  $\frac{1}{2}\int_{\Omega_0} (u^\star)^2 dx$  & $\PP(F(u^\star,\xi)\geq z)$   & Time [sec]\\ \midrule
                \multicolumn{4}{c}{\textit{LDT} $(k = 1)$} \\
				$10^{-1}$ & $3.71$ &  $4.9\cdot 10^{-2}$  & $14.8$ \\
				$10^{-2}$ & $4.68$ &  $5.4\cdot 10^{-3}$  & $15.1$ \\
				$10^{-3}$ & $5.49$ &  $3.9\cdot 10^{-4}$  & $14.4$ \\
				$10^{-4}$ & $6.12$ &  $4.5\cdot 10^{-5}$  & $17.5$ \\ 
				$10^{-5}$ & $6.63$ &  $2.8\cdot 10^{-6}$  & $21.8$ \\
				$10^{-6}$ & $6.97$ &  $3.0\cdot 10^{-7}$  & $31.7$ \\ 
                \midrule
                \multicolumn{4}{c}{\textit{LDT} $(k = 2)$} \\
                $10^{-1}$ & $3.52$ &  $8.9\cdot 10^{-2}$  & $3.5 \cdot 10^{5}$ \\
                $10^{-2}$ & $4.43$ &  $9.5\cdot 10^{-3}$  & $3,430.0$ \\
                $10^{-3}$ & $5.22$ &  $9.2\cdot 10^{-4}$  & $3,899.0$ \\
                $10^{-4}$ & $5.88$ &  $9.5\cdot 10^{-5}$  & $3,913.0$ \\ 
                $10^{-5}$ & $6.41$ &  $9.6\cdot 10^{-6}$  & $4,137.5$ \\
                $10^{-6}$ & $6.80$ &  $9.4\cdot 10^{-7}$  & $5,069.5$ \\ 
				\bottomrule 
			\end{tabular} 
		\end{center}
	}
\end{table}

\begin{figure}[tbhp]\centering
	\begin{tikzpicture}[]
	\begin{axis}[width=0.6\textwidth, 
	height = 0.3\textwidth,
	scale only axis, xmin=0, xmax=1,
legend pos=outer north east,
	legend style ={font =\small,
	}, 
	xlabel=$x_2$,
	ylabel = $u^\star$,
	cycle list name= 
	my mark
	,
	]
	\addplot coordinates {
		(0.000e+00, -2.334e+00)  (3.448e-02, -1.191e+00)  (6.897e-02, -1.239e+00)  (1.034e-01, -1.309e+00)  (1.379e-01, -1.398e+00)  (1.724e-01, -1.504e+00)  (2.069e-01, -1.619e+00)  (2.414e-01, -1.738e+00)  (2.759e-01, -1.851e+00)  (3.103e-01, -1.948e+00)  (3.448e-01, -2.016e+00)  (3.793e-01, -2.042e+00)  (4.138e-01, -2.014e+00)  (4.483e-01, -1.920e+00)  (4.828e-01, -1.753e+00)  (5.172e-01, -1.507e+00)  (5.517e-01, -1.177e+00)  (5.862e-01, -7.573e-01)  (6.207e-01, -4.588e+00)  (6.552e-01, -4.253e+00)  (6.897e-01, -3.983e+00)  (7.241e-01, -3.756e+00)  (7.586e-01, -3.560e+00)  (7.931e-01, -3.391e+00)  (8.276e-01, -3.247e+00)  (8.621e-01, -3.127e+00)  (8.966e-01, -3.031e+00)  (9.310e-01, -2.959e+00)  (9.655e-01, -2.911e+00)  (1.000e+00, -5.773e+00) 
	};
	\addlegendentry{$\alpha=10^{-1}
		$}
	\addplot coordinates {
		 (0.000e+00, -1.823e+00)  (3.448e-02, -9.526e-01)  (6.897e-02, -1.034e+00)  (1.034e-01, -1.155e+00)  (1.379e-01, -1.313e+00)  (1.724e-01, -1.503e+00)  (2.069e-01, -1.719e+00)  (2.414e-01, -1.950e+00)  (2.759e-01, -2.179e+00)  (3.103e-01, -2.388e+00)  (3.448e-01, -2.552e+00)  (3.793e-01, -2.644e+00)  (4.138e-01, -2.639e+00)  (4.483e-01, -2.513e+00)  (4.828e-01, -2.252e+00)  (5.172e-01, -1.848e+00)  (5.517e-01, -1.298e+00)  (5.862e-01, -6.031e-01)  (6.207e-01, -5.281e+00)  (6.552e-01, -4.820e+00)  (6.897e-01, -4.466e+00)  (7.241e-01, -4.180e+00)  (7.586e-01, -3.939e+00)  (7.931e-01, -3.734e+00)  (8.276e-01, -3.562e+00)  (8.621e-01, -3.420e+00)  (8.966e-01, -3.307e+00)  (9.310e-01, -3.223e+00)  (9.655e-01, -3.167e+00)  (1.000e+00, -6.278e+00) 
	};
	\addlegendentry{$\alpha=10^{-2}
		$}
	\addplot coordinates {
(0.000e+00, -1.258e+00)  (3.448e-02, -6.807e-01)  (6.897e-02, -7.853e-01)  (1.034e-01, -9.440e-01)  (1.379e-01, -1.158e+00)  (1.724e-01, -1.425e+00)  (2.069e-01, -1.741e+00)  (2.414e-01, -2.093e+00)  (2.759e-01, -2.463e+00)  (3.103e-01, -2.819e+00)  (3.448e-01, -3.123e+00)  (3.793e-01, -3.328e+00)  (4.138e-01, -3.387e+00)  (4.483e-01, -3.257e+00)  (4.828e-01, -2.909e+00)  (5.172e-01, -2.332e+00)  (5.517e-01, -1.529e+00)  (5.862e-01, -5.171e-01)  (6.207e-01, -5.730e+00)  (6.552e-01, -5.163e+00)  (6.897e-01, -4.744e+00)  (7.241e-01, -4.414e+00)  (7.586e-01, -4.142e+00)  (7.931e-01, -3.914e+00)  (8.276e-01, -3.724e+00)  (8.621e-01, -3.568e+00)  (8.966e-01, -3.445e+00)  (9.310e-01, -3.353e+00)  (9.655e-01, -3.292e+00)  (1.000e+00, -6.524e+00) 
	};
	\addlegendentry{$\alpha=10^{-3}
		$}
	\addplot coordinates {
(0.000e+00, -7.733e-01)  (3.448e-02, -4.395e-01)  (6.897e-02, -5.484e-01)  (1.034e-01, -7.192e-01)  (1.379e-01, -9.589e-01)  (1.724e-01, -1.274e+00)  (2.069e-01, -1.665e+00)  (2.414e-01, -2.126e+00)  (2.759e-01, -2.636e+00)  (3.103e-01, -3.159e+00)  (3.448e-01, -3.638e+00)  (3.793e-01, -4.002e+00)  (4.138e-01, -4.172e+00)  (4.483e-01, -4.075e+00)  (4.828e-01, -3.659e+00)  (5.172e-01, -2.904e+00)  (5.517e-01, -1.822e+00)  (5.862e-01, -4.550e-01)  (6.207e-01, -5.977e+00)  (6.552e-01, -5.323e+00)  (6.897e-01, -4.857e+00)  (7.241e-01, -4.498e+00)  (7.586e-01, -4.206e+00)  (7.931e-01, -3.965e+00)  (8.276e-01, -3.765e+00)  (8.621e-01, -3.602e+00)  (8.966e-01, -3.473e+00)  (9.310e-01, -3.378e+00)  (9.655e-01, -3.315e+00)  (1.000e+00, -6.568e+00) 
	};
	\addlegendentry{$\alpha=10^{-4}
		$}
	\addplot coordinates {
(0.000e+00, -4.295e-01)  (3.448e-02, -2.602e-01)  (6.897e-02, -3.565e-01)  (1.034e-01, -5.144e-01)  (1.379e-01, -7.482e-01)  (1.724e-01, -1.074e+00)  (2.069e-01, -1.504e+00)  (2.414e-01, -2.042e+00)  (2.759e-01, -2.676e+00)  (3.103e-01, -3.365e+00)  (3.448e-01, -4.040e+00)  (3.793e-01, -4.600e+00)  (4.138e-01, -4.927e+00)  (4.483e-01, -4.907e+00)  (4.828e-01, -4.451e+00)  (5.172e-01, -3.526e+00)  (5.517e-01, -2.152e+00)  (5.862e-01, -4.047e-01)  (6.207e-01, -6.051e+00)  (6.552e-01, -5.334e+00)  (6.897e-01, -4.838e+00)  (7.241e-01, -4.463e+00)  (7.586e-01, -4.162e+00)  (7.931e-01, -3.916e+00)  (8.276e-01, -3.712e+00)  (8.621e-01, -3.548e+00)  (8.966e-01, -3.419e+00)  (9.310e-01, -3.323e+00)  (9.655e-01, -3.260e+00)  (1.000e+00, -6.457e+00)  
	};
	\addlegendentry{$\alpha=10^{-5}
		$}
	\addplot coordinates {
		(0.000e+00, -2.190e-01)  (3.448e-02, -1.435e-01)  (6.897e-02, -2.184e-01)  (1.034e-01, -3.481e-01)  (1.379e-01, -5.528e-01)  (1.724e-01, -8.581e-01)  (2.069e-01, -1.290e+00)  (2.414e-01, -1.867e+00)  (2.759e-01, -2.591e+00)  (3.103e-01, -3.427e+00)  (3.448e-01, -4.299e+00)  (3.793e-01, -5.077e+00)  (4.138e-01, -5.600e+00)  (4.483e-01, -5.697e+00)  (4.828e-01, -5.237e+00)  (5.172e-01, -4.161e+00)  (5.517e-01, -2.497e+00)  (5.862e-01, -3.618e-01)  (6.207e-01, -5.985e+00)  (6.552e-01, -5.227e+00)  (6.897e-01, -4.717e+00)  (7.241e-01, -4.338e+00)  (7.586e-01, -4.038e+00)  (7.931e-01, -3.793e+00)  (8.276e-01, -3.592e+00)  (8.621e-01, -3.430e+00)  (8.966e-01, -3.304e+00)  (9.310e-01, -3.210e+00)  (9.655e-01, -3.148e+00)  (1.000e+00, -6.235e+00)  
	};
	\addlegendentry{$\alpha=10^{-6}
		$}
	\end{axis}
	\end{tikzpicture}
	\caption{The optimal boundary condition $u^\star$ from problem \eqref{eq:sad-minreg} with $z=1$ for different $\alpha$'s using the second-order LDT approximation.}\label{fig:sad-reg30-u}
\end{figure}
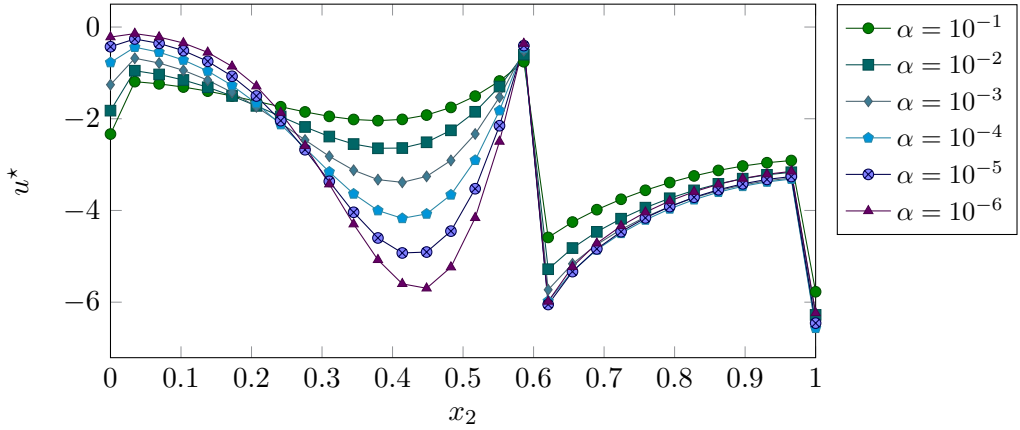

\begin{figure}[tbhp]\centering
	\begin{tikzpicture}[]
	\begin{groupplot}[group style = {columns=3, 
			rows=1, 
		vertical sep =50pt, horizontal sep = 30pt}, width = 5cm, height = 4.0cm,view/h=300, point meta max=4, point meta min=-6, zmin=-6.4,zmax=5,xlabel=$x_1$, ylabel=$x_2$,	
	tick label style={font=\footnotesize},  title style = {font= \small},
	label style={font=\footnotesize},
	]
	\nextgroupplot[
	title={ $\alpha=10^{-2}
		$},]
	\addplot3[surf, mesh/rows=30] table[skip first n=0, z expr=\thisrow{Z}]{\data/reg_n30z1alpha0.01.txt};
	\nextgroupplot[
	title={ $\alpha=10^{-4}
		$},]
	\addplot3[surf, mesh/rows=30] table[skip first n=0, z expr=\thisrow{Z}]{\data/reg_n30z1alpha0.0001.txt};
	\nextgroupplot[colorbar right,
	title={ $\alpha=10^{-6}
		$},]
	\addplot3[surf, mesh/rows=30] table[skip first n=0, z expr=\thisrow{Z}]{\data/reg_n30z1alpha1.0e-6.txt};
	\end{groupplot}
	\end{tikzpicture}
	\caption{The temperature profile $y$ under the optimal boundary condition $u^\star$ in \Cref{fig:sad-reg30-u} and $\xistar$ for the minimization problem \eqref{eq:sad-minreg} with $z=1$ for different $\alpha$'s using the second-order LDT approximation.
	}\label{fig:sad-reg30-y}
\end{figure}
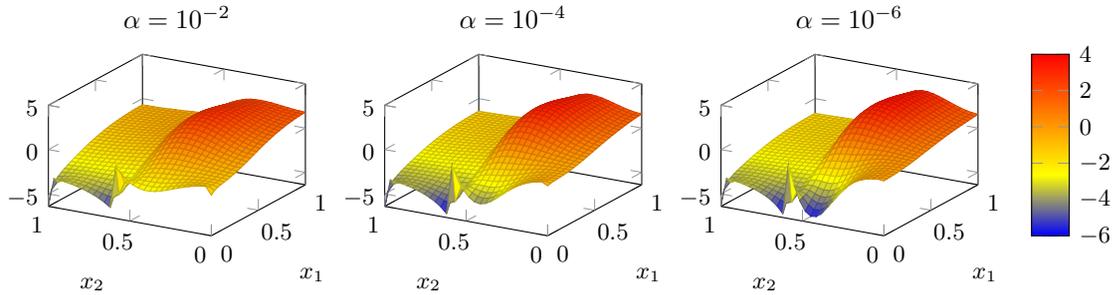

We use a finite-difference method to discretize the PDE \eqref{eq:sad-PDE} with $m=30$ elements.
Therefore, the mesh is $30\times30$ and the state $y$ is 900-dimensional. %
To avoid issues of infeasibility and convergence to non-optimal stationary points, we successively warm-start the solution of the $\alpha=10^{-(k+1)}$ problem using the optimal values determined for the $\alpha=10^{-k}$ problem.
In case of the first-order LDT approximation, we use a simple finite difference method to compute gradients, since the random parameter is low-dimensional $n=2$ (this also explains the faster solution times in \Cref{tab:sad1}).
In case of the second-order LDT approximation, we use forward-mode automatic differentiation to compute gradients, since we need accurate gradients of the correction term \eqref{eq:op-G-SO} to ensure convergence.

\Cref{tab:sad1} summarizes the computational performance.
We observe that the optimal objective values become larger for more extreme events; \textit{i.e.}, the optimal boundary control $u^\star$ has larger variations.
The terms $\PP(F(u^\star,\xi)\geq z)$ show chance constraint feasibilities, and the fact that they are all less than $\alpha$ indicates that all obtained solutions are feasible.  %
Comparing the two approaches, we observe that the optimal objective values of the second-order approximations are better compared with those of first-order approximations;
however, this improvement comes at a significantly higher computational cost.
We note that a large portion of the computational time is due to the inability of forward-mode automatic differentiation in JuMP to differentiate a sparse solve that required us to modify the PDE solve to a dense solve.
The latter can be reduced by using a more efficient automatic differentiation implementation (or perhaps by using tailored codes for evaluating adjoints and gradients).

Finally, \Cref{fig:sad-reg30-u,fig:sad-reg30-y} show the optimal boundary control $u^\star$ and corresponding temperature profile $y$ under the optimal $u$ and $\xistar$ obtained through \eqref{eq:op-G}.
The optimal control inputs $u$ have similar shapes but scale with $\log \alpha$; moreover, they all have a jump at $x_2=0.6$, which comes from the discontinuity of the diffusion coefficient $\kappa$ at $x_2=0.6$.
The temperature profile $y$ under the optimal boundary condition $u^\star$ and the LDT realization $\xistar$, corresponds to the most representative profile among all feasible profiles in \eqref{eq:sad-minreg}.
Similar to the boundary conditions in \Cref{fig:sad-reg30-u}, the temperatures for smaller risk thresholds have larger variations.

%% file: conclusions.tex
\section{Discussions and conclusions}\label{sec:conclusions}

Despite their rare occurrence, factors such as climate change and human population growth are causing a steady increase in the frequency of extreme high-impact events. The paucity of available data and inability to accurately model socioeconomic costs during these events motivates chance constraints as a viable modeling framework to mitigate the risk of extreme events.

In this paper, we combined ideas from large deviation theory and convex and bilevel optimization to propose a novel solution method for decision-making problems constrained by probabilities of rare events.
Unlike classical approaches, our formulations are sampling-free, independent of event extremeness, and applicable to a broad class of nonlinear problems.
Under certain regularity assumptions on the system constraints,
they constitute safe conservative approximations or even asymptotically equivalent reformulations of the true problem.
Furthermore, by utilizing local Taylor approximations of the system functions, we provide refined approximations over classical large deviation estimates that turn out to be empirically accurate even in non-asymptotic regimes.
Although we only considered problems affected by Gaussian mixture uncertainties, it should be noted that our methods are also directly applicable to more general settings where the uncertainties can be transformed to a Gaussian or approximated well using a mixture of Gaussians, for which there are efficient off-the-shelf codes.
Our proposed formulations complement this well since they can also be solved by off-the-shelf optimization solvers.
Computational experiments on applications from diverse domains confirm the broad applicability, and improved efficiency and accuracy of our method over classical sampling-based methods in the rare event regime:
the portfolio application illustrates its applicability to nonlinear high-dimensional uncertainties, the structural design application showcases its potential advantages over sampling-based methods in systems that are nonlinear in the decisions as well as uncertainties, whereas the optimal control application shows its ability to address complicated PDE-constrained systems. %

We believe that our approach is but a small step towards tackling rare chance constraints, and much more work needs to be done. Various regularity assumptions on the structure of the constraint function (\textit{e.g.}, smoothness, concavity etc.) as well as the uncertainty (\textit{e.g.}, sufficiently smooth rate functions) need to be lifted and generalized. For example, our findings suggest that we cannot expect our method, and perhaps even large deviation theory, to work well when the system constraints are convex functions of the uncertainty. %
We would also like to extend our methodology to nonsmooth constraint functions, where the current Taylor-based probability estimates would be inapplicable, as well as to joint chance constrained problems. %
From a computational viewpoint, although our method has the advantage of using off-the-shelf codes, this also means that it can be inefficient at tackling more structured large-scale problems.
For example, we need more efficient methods to compute derivatives of the probability estimates, and particularly of the second-order correction term. We believe this may also spur the development of new techniques in automatic differentiation.